\documentclass[moor]{informs1}              

\usepackage{natbib}
 \NatBibNumeric
 \def\bibfont{\small}%
 \bibpunct[, ]{[}{]}{,}{n}{}{,}%
 
 \usepackage{graphicx}
 \usepackage{url,color}
 \usepackage{subfigure}
 \usepackage{amsfonts,mathrsfs}
 \usepackage{amssymb,amsmath}
 \usepackage{verbatim}
 \usepackage{acronym}
 \usepackage{multirow}
 \usepackage{array}
 \usepackage[utf8]{inputenc}
 \usepackage[english]{babel}
 \usepackage[mathscr]{euscript}
 \usepackage{bbold}
 \usepackage{bbm}

\usepackage{textcomp}%
\usepackage{manyfoot}%
\usepackage{booktabs}%
\usepackage{algorithm}%
\usepackage{algorithmicx}%
\usepackage{algpseudocode}%
\usepackage{listings}%

 \newcommand{\red}{\textcolor{black}}
 
\usepackage[colorlinks=true,breaklinks=true,bookmarks=true,urlcolor=blue,
     citecolor=blue,linkcolor=blue,bookmarksopen=false,draft=false]{hyperref}

\def\EMAIL#1{\href{mailto:#1}{#1}}

\TheoremsNumberedThrough     

\EquationsNumberedThrough    

\begin{document}

\TITLE{A stochastic perturbed augmented Lagrangian method for smooth convex constrained minimization}

\ARTICLEAUTHORS{%
\AUTHOR{Nitesh Kumar Singh}
\AFF{School of Computer Science, University of Petroleum and Energy Studies, Dehradun, Uttarakhand, 248007 India
    \EMAIL{niteshkumar.singh@ddn.upes.ac.in.} 
    }
\AUTHOR{Ion  Necoara}
\AFF{Automatic Control and Systems
	Engineering Department, University Politehnica Bucharest, 060042
	Bucharest, Romania  and  Gheorghe Mihoc-Caius Iacob  Institute of Mathematical Statistics and Applied Mathematics of the Romanian Academy, 050711 Bucharest, Romania, \EMAIL{ion.necoara@upb.ro.} 
    }
} 

\ABSTRACT{%
This paper considers  smooth convex  optimization problems with  many functional constraints. To solve this general  class of problems we propose  a new stochastic  perturbed augmented Lagrangian  method, called SGDPA, where a perturbation is introduced in the augmented Lagrangian function by multiplying the dual variables with a subunitary parameter. Essentially, we linearize the objective and one randomly chosen functional constraint within the perturbed augmented Lagrangian at the current iterate and add a quadratic regularization that leads to a stochastic gradient descent update for the primal variables, followed by a perturbed random coordinate  ascent step to update the dual variables.   We provide a convergence analysis in both optimality and feasibility criteria for the iterates of SGDPA algorithm using basic assumptions on the problem. In particular,  when the dual updates are assumed to be bounded, we prove sublinear rates of convergence for the iterates of  algorithm SGDPA  of order $\mathcal{O} (k^{-1/2})$ when the objective is convex and of order $\mathcal{O} (k^{-1})$ when the objective is strongly convex, where $k$ is the iteration counter. Under some additional assumptions, we  prove that the dual iterates are bounded and in this case we  obtain  convergence rates of order $\mathcal{O} (k^{-1/4})$  and  $\mathcal{O} (k^{-1/2})$ when the objective is convex and strongly convex, respectively. Preliminary numerical experiments on  problems with many quadratic constraints demonstrate the viability and performance of our method when compared to some existing state-of-the-art optimization methods and software.
}

\MSCCLASS{90C25, 90C06, 65K05.}
\HISTORY{December 2024}

\maketitle

\section{Introduction}
This paper addresses the intricate challenges associated with the minimization of smooth convex objective functions and complex feasible sets. The central focus is on the following convex constrained optimization problem:
\begin{equation} \label{eq:prob}
\begin{array}{rl}
F^* = & \min\limits_{x \in \mathcal{Y} \subseteq \mathbb{R}^n} \; F(x) \\
& \text{subject to } \;  h_j(x) \le 0 \;\; \forall j = 1:m,
\end{array}
\end{equation}
where the objective function $F$ and the functional constraints $h_j$, for all $j = 1:m $, are assumed to be convex and smooth.  Additionally, $\mathcal{Y}$ is a nonempty, closed, convex and simple set (by simple we mean that the projection onto this set is easy).  Such problems pose significant computational and theoretical challenges, particularly when dealing with \textit{large} number of constraints (i.e., $m$~large).

\medskip 

\noindent \textit{Motivation}. 
Problems of the form \eqref{eq:prob} appear, e.g., in robust optimization. In contrast to the chance-constrained problems, which considers probabilistic constraints, robust optimization requires that the constraint \( h(\mathbf{x}; \xi) \leq 0 \) be satisfied for every scenario \(\xi \in \Omega\), where $\Omega$ can be even an infinite set. Hence, robust optimization problems are  formulated as follows:
\begin{align} 
\label{eq:robust_opt}
\min_{\mathbf{x} \in \mathcal{Y} \subseteq \mathbb{R}^n} \, F(\mathbf{x}), \quad \text{s.t.} \quad h(\mathbf{x}, \xi) \leq 0 \quad \forall \xi \in \Omega.
\end{align}

\noindent Similar to the scenario approximation technique used for chance-constrained problems (see, e.g., \cite{CalCam:05}), a sampling approach has been proposed to numerically solve this robust optimization problem. Suppose \(\{ \xi_1, \ldots, \xi_m \}\) are \(m\) independently drawn samples and consider the following optimization problem with a finite number of constraints:
\begin{align} 
\label{eq:robust_opt_m}
\min_{\mathbf{x} \in \mathcal{Y} \subseteq \mathbb{R}^n} \, F(\mathbf{x}), \quad \text{s.t.} \quad h(\mathbf{x}, \xi_i) \leq 0 \quad \forall i=1:m.
\end{align}
It has been proved in  \cite{CalCam:05} that, for any \(\delta \in (0, 1)\) and any \(\varepsilon \in (0, 1)\), if the number of samples satisfies \(m \geq \frac{n}{\delta \varepsilon}\), then any  solution of problem \eqref{eq:robust_opt_m}  is a \(\delta\)-level robustly feasible solution  for original problem  \eqref{eq:robust_opt} with probability at least \(1 - \varepsilon\). However, if \(n\) is large and high feasibility levels along with high probabilities are required, \(m\) will be very large, leading to an optimization problem of the form \eqref{eq:prob}  with a  large number of functional constraints. Similarly, model predictive control problems for linear dynamical systems with state and input constraints lead to the optimization problem \eqref{eq:prob}  with a  large number of functional constraints, see  \cite{JiaStu:05, NedNec:14}. 

\medskip

\noindent \textit{Related work}. 
In the field of optimization, various  \textit{primal-dual} (i.e., updating both, primal and dual variables) \textit{deterministic} (i.e., working with the full set of constraints) methods   have been developed to solve smooth optimization problems of the form \eqref{eq:prob}. A notable contribution is the work \cite{Xu:21}, which addresses smooth convex optimization problems involving a composite objective function along with linear equality and nonlinear inequality constraints.  Based on the standard augmented Lagrangian, \cite{Xu:21} proposed the Linearized Augmented Lagrangian Method (LALM). In each iteration of LALM, a single proximal gradient step is performed for updating the primal variables, followed by a gradient ascent update for the dual variables. It has demonstrated in \cite{Xu:21} that LALM achieves a sublinear rate of convergence of order $\mathcal{O} (k^{-1})$, where $k$ is the iteration counter. For convex problems with nonlinear constraints, papers  \cite{LiuMa:19,NedNec:14} examined  inexact augmented Lagrangian algorithms, utilizing Nesterov’s optimal first-order method \cite{Nes:13} to approximately solve each primal subproblem. These papers also provide accuracy requirements for solving each primal subproblem to estimate the outer iteration complexity of the corresponding inexact augmented Lagrangian algorithm. In particular,  it is shown that $\mathcal{O} (\epsilon^{-1})$ gradient evaluations are sufficient to obtain a primal $\epsilon$ solution (see definition  below).  

\medskip 

\noindent Other works focused on the perturbed augmented Lagrangian framework, proposing  deterministic perturbed augmented Lagrangian methods. A key algorithm that inspired also our current work is the Gradient Descent and Perturbed Ascent (GDPA) method introduced in \cite{Lu:22}.  Specifically tailored for smooth constrained problems, GDPA is a single-loop algorithm that can find KKT points of nonconvex optimization problems under nonconvex inequality constraints. At each iteration of GDPA,  the linearized perturbed augmented Lagrangian function with a proximal term is minimized to update the primal variables, followed by a perturbed full dual update. In this update, the perturbation technique is introduced by adding a negative curvature to the maximization problem, which ensures that the dual update is well-behaved. GDPA achieves a convergence rate of order $\mathcal{O}\left(k^{-1/3}\right)$ to KKT points. Several other studies have also considered this framework, designed algorithms for solving optimization problems with nonconvex objective under linear equality constraints or convex constraints, as seen in \cite{KosSha:11, KopRib:17, YanSpa:22, HajHon:19, HonRaz:16, ZhaLuo:20}. The perturbed dual step, can be interpreted as performing a dual ascent on a certain regularized Lagrangian in the dual space, see   \cite{KopRib:17}. The primary purpose of introducing dual perturbation/regularization in this context, and in many related works, is to ensure that the dual updates remain well-behaved and  easier to analyze. Intuitively, adopting and modifying this perturbation strategy  eliminates the need to bound the difference of successive dual variables relative to the difference of successive primal variables, as the change in dual variables is now well controlled (see \cite{HajHon:19}). Conversely, without perturbation, bounding the difference of successive dual variables is not straightforward (see Lemma 1 in \cite{HonZha:17}). Moreover, as demonstrated in \cite{YanSpa:22}, this dual perturbation technique allows for proving the boundedness of dual updates, leading to the establishment of a sufficiently decreasing and lower-bounded Lyapunov function, which is crucial for establishing convergence. 
Although the algorithms developed in these works achieve the standard convergence rate for a given class of problems, they still face a significant technical challenge: they must consider the \textit{entire} set of constraints in the subproblem, rendering them intractable when the number of constraints is large.

\medskip 

\noindent It is known that \textit{stochastic} methods can overcome such technical problem. Stochastic (sub)gradient framework is a prevalent methodology for minimizing finite sum objective functions  \cite{NemYud:83}. When combined with simple constraints, and the computation of the prox or projection operators is straightforward, a plethora of methods emerge, such as stochastic gradient descent (SGD) and stochastic proximal point (SPP) algorithms \cite{MouBac:11, NemJud:09, Nec:20, WanBer:16}.  Numerous extensions of SGD/SPP have been devised for solving convex problem \eqref{eq:prob} with \textit{nonsmooth} functional constraints, see e.g., \cite{Ned:11, NecSin:22}, and their mini-batch variants \cite{SinNec:23, NecNed:21}. Additionally,  a Lagrangian primal-dual stochastic subgradient (PDSG) method for convex nonsmooth functional constraint problems have been designed in \cite{Xu:20}. The results from \cite{Xu:20} can also be extended to the smooth constrained case using only tools from nonsmooth optimization. However,  it is important to note that  PDSG algorithm uses the previous primal iterate to update the dual variables, instead of considering the latest primal update, thus making the algorithm to be slow. All these algorithms achieve sublinear rates of convergence in terms of optimality and feasibility for an average sequence when they are combined with decreasing stepsize rules. However, extensions of stochastic first-order methods to problems with \textit{smooth} functional constraints are limited. For example, a single-loop stochastic augmented Lagrangian primal-dual method is proposed in \cite{JinWan:22} for solving nonconvex problem \eqref{eq:prob} with smooth functional constraints. This method achieves a convergence rate of order  $\mathcal{O} (k^{-1/4})$ for finding an approximate KKT point. In cases where the initial point is feasible, this  reduces to $\mathcal{O} (k^{-1/3})$. However, to ensure the boundedness of dual multipliers, the penalty parameter in augmented Lagrangian and the stepsizes must be chosen dependent on the final iteration counter, hence it needs to be fixed a priori. In this paper, based on the perturbed augmented Lagrangian framework, we focus on developing a stochastic linearized perturbed augmented Lagrangian algorithm that considers only one randomly selected constraint at each iteration and leverages the dual perturbation technique, particularly for bounding the dual variables.

\medskip

\noindent \textit{Contributions}.  
This paper introduces a novel primal-dual method for solving smooth convex  problems with many functional  constraints, named \textit{Stochastic Gradient Descent Perturbed Ascent} (SGDPA) algorithm. To the best of our knowledge, this work is one of the first proposing a pure stochastic gradient descent and coordinate perturbed ascent  which is supported by  efficiency estimates for smooth convex constrained optimization problems.   The main  contributions are:

(i) We consider a general smooth convex optimization model that has a large number of functional constraints. For solving this problem,  we consider a perturbed augmented Lagrangian approach  and  we introduce a new primal-dual algorithm, called  SGDPA, that combines two approaches:  a \textit{stochastic gradient} descent step w.r.t. the perturbed augmented Lagrangian based only on a randomly chosen constraint, and then  a  perturbed dual \textit{random coordinate gradient} ascent step to update the dual variables. For SGDPA in both steps we choose the stochastic indexes uniformly at random and independently. To perform the dual step we  choose a fixed positive penalty parameter and a fixed subunitary perturbation parameter. Being stochastic in nature, the computational simplicity of the SGDPA subproblem, which uses only gradient information and considers one constraint at a time, makes it suitable for problems with many functional constraints. Additionally, SGDPA's ability to operate without requiring an initial feasible point is an advantage over existing algorithms such as  \cite{AusTeb:10,JinWan:22,Lu:22}.

(ii) It is well-known that proving boundedness of the multiplier sequence in Lagrangian based schemes is a difficult matter. In this paper, depending on the boundedness nature of the dual updates, i.e., whether they are assumed to be bounded or, under some assumptions, proven to be bounded, we provide a detailed convergence analysis for SGDPA, yielding  convergence rates in both optimality and feasibility criteria evaluated at some average point. Specifically, when the dual updates are assumed to be bounded, we derive \textit{sublinear convergence rates} of order $\mathcal{O} (k^{-1/2})$ when the objective is convex and of order $\mathcal{O} (k^{-1})$ when the objective is strongly convex. Under some additional assumptions, we  prove that the dual iterates are bounded, but the bound depends on the perturbation parameter,  and in this case we  obtain  convergence rates of order $\mathcal{O} (k^{-1/4})$  and of order $\mathcal{O} (k^{-1/2})$ when the objective is convex and strongly convex, respectively.  To the best of our knowledge, these results are novel and align with the optimal convergence rates for stochastic first-order methods  reported in the literature.  

(iii) Finally, detailed numerical experiments on  problems with many quadratic constraints using synthetic and real data demonstrate the viability and performance of our method when compared to some existing state-of-the-art optimization methods and software.

\medskip 

\noindent \textit{Content}.  
The remainder of this paper is organized as follows: In Section 2, we introduce the essential notations and the key assumptions. Section 3, presents the perturbed augmented Lagrangian formulation and its properties. Section 4 introduces the new algorithm. In Section 5, we delve into the convergence analysis of the proposed  algorithm. Finally, in Section 6, we present numerical results that validate the practical efficacy of our algorithm.


\section{Notations and assumptions}
Let us define the individual sets $\mathcal{X}_j=\left\{ x \in  \mathbb{R}^n | \;  h_j(x)\le 0 \right\}$, for all $j = 1:m$. Then, the feasible set of the optimization problem \eqref{eq:prob} can be written as $\mathcal{X}= \mathcal{Y} \cap \left(\cap_{j = 1:m} \mathcal{X}_j \right)$.   We assume that the feasible set $\mathcal{X}$ is nonempty and  problem \eqref{eq:prob} has a finite optimum, i.e., the problem is well-posed. Let  $F^*$ and $\mathcal{X}^*$ denote the optimal value and the optimal set of \eqref{eq:prob}, respectively:
\[F^*=\min_{x  \in \mathcal{X}} F(x), \quad \mathcal{X}^*=\{x\in \mathcal{X} \mid F(x)=F^*\} \not= \emptyset.\]
Moreover, $x^*$ denotes a primal solution of the primal problem  \eqref{eq:prob}, i.e., $x^* \in \mathcal{X}^*$,  and $\lambda^* \in \mathbb{R}_+^m$ denotes a  dual solution (i.e., an optimal point of the corresponding dual problem). To denote the gradient of any function w.r.t. the primal variable $x$, we use $\nabla$ instead of $\nabla_x$. For any given real $a$, we also use the notation $(a)_+ = \max(a,0)$. We denote $\mathcal{N}_\mathcal{Y} (x)$  the normal cone of the nonempty closed convex set  $\mathcal{Y}$ at $x$ and $\Pi_\mathcal{Y} (x)$ the projection of $x$ onto the set $\mathcal{Y}$.



\medskip 

\noindent Now, let us introduce our main assumptions. First, we assume that the following Karush-Kuhn-Tucker (KKT) conditions hold for problem \eqref{eq:prob}:
\begin{assumption}\label{Ass:KKT}
 For any $\tau\ge0$, there exists a primal-dual solution $(x^*, \lambda^*)$ for \eqref{eq:prob}  satisfying the Karush-Kuhn-Tucker (KKT) conditions:
    \begin{align}
        & 0 \in \nabla F(x^*) + \mathcal{N}_\mathcal{Y} (x^*) + \frac{1-\tau}{m} \sum_{j = 1}^m (\lambda^*)^j \nabla h_j(x^*), \label{eq:KKT1}\\
        & x^* \in \mathcal{Y}, \;\; h_j(x^*) \le 0, \;\;  (\lambda^*)^j \ge 0, \;\; (\lambda^*)^j h_j(x^*) =  0 \;\;\;  \forall j = 1:m. \label{eq:KKT3}
    \end{align}
\end{assumption}

\noindent Assumption \ref{Ass:KKT} is satisfied if e.g.,  a certain constraint qualification condition holds, such as the Slater's condition,  \cite{KosSha:11, Xu:20}. Further, we also assume  smoothness for the objective function  and the functional constraints.
\begin{assumption}\label{Ass:Smooth}
    The objective function $F$ and the functional constraints $h_j$, for all $j = 1:m$, are convex and continuously differentiable on the convex set $\mathcal{Y}$ with $L_f>0$ and $L_h>0$ Lipschitz continuous gradients, respectively, i.e.,  we have:
    \begin{align*}
        & \|\nabla F(y) - \nabla F(x) \| \le {L_f} \|y-x\| \quad \;\; \forall x,y \in \mathcal{Y}, \\
        & \|\nabla h_j(y) - \nabla h_j(x) \| \le {L_h} \|y-x\|\quad \forall x,y \in \mathcal{Y}, \;\; j=1:m.
    \end{align*}
\end{assumption}

\noindent In our convergence analysis  we also use  the following bounds:  there exists $B_F, M_h, B_h > 0$ such that
\begin{align}\label{eq:Bddness}
    \|\nabla F(x)\|\le B_F, \quad |h_j(x)| \le M_h,\quad \|\nabla h_j(x)\| \le B_h \quad  \forall x \in \mathcal{Y}, \; j = 1:m.
\end{align}
We use these bounds, \eqref{eq:Bddness},  to prove the smoothness property of the perturbed augmented Lagrangian \eqref{eq:AugLag} and also  the boundedness of the dual iterates (see Section \ref{ch04:S2sss} below). Note that \eqref{eq:Bddness} always holds if e.g., $\mathcal{Y}$ is a bounded set, since the functions $F$ and $h_j$'s are assumed  continuously differentiable.  
Our next assumption is related to the strong convexity condition of the objective $F$:
\begin{assumption}
\label{Ass:StrongConv}
The objective function $F$ is $\mu$-strongly convex, i.e., there exists  $\mu \ge 0$ such that the following relation holds:
   \begin{align*}
       F(y) \ge F(x) + \langle \nabla F(x), y-x \rangle + \frac{\mu}{2} \|y-x\|^2 \;\;\;\; \forall x,y \in \mathcal{Y}.
   \end{align*}
\end{assumption}
\noindent Note that if $\mu = 0$, then $F$ is simply a convex function.
Let us also define a stochastic approximately optimal solution for problem \eqref{eq:prob}.

\begin{definition}
A random vector $x \in \mathcal{Y}$ is called a stochastic $\epsilon$-approximately optimal solution of problem \eqref{eq:prob} if:
    $$\mathbb{E} [|F(x) - F(x^*)|]\le \epsilon \;\; \text{ and } \;\; \mathbb{E} \left[ \frac{1}{m} \sum_{j=1}^m [h_j(x)]_+ \right] \le \epsilon.$$
\end{definition}


\section{Perturbed augmented Lagrangian and its properties}
\noindent Let us now introduce the perturbed augmented Lagrangian function for the optimization problem \eqref{eq:prob}, inspired from  \cite{Lu:22}. This function is expressed as follows: 
\begin{align}\label{eq:AugLag}
    \mathcal{L}_{\rho,\tau} (x;\lambda) & := 
    F(x) + \Psi_{\rho,\tau} (h(x);\lambda),
\end{align}
where the term $\Psi_{\rho,\tau} (h(x);\lambda)$ is defined as
$$ \Psi_{\rho,\tau} (h(x);\lambda) = \frac{1}{m}\sum_{j = 1}^m \psi^j_{\rho,\tau} (h_j(x); \lambda^j), $$  
\noindent with each component $\psi^j_{\rho,\tau} (h_j(x); \lambda^j)$ expressed as
$$ \psi^j_{\rho,\tau} (h_j(x); \lambda^j) = \frac{1}{2 \rho} \left[ (\rho h_j(x) + (1 - \tau) \lambda^j)^2_+ - ((1 - \tau) \lambda^j)^2 \right]. $$
\noindent Here, $\rho > 0$ is the penalty parameter, $\tau\in [0,1)$ is the perturbation parameter, and $\lambda$ is the Lagrangian multiplier or dual variable. The function $\mathcal{L}_{\rho,\tau} (x;\lambda)$ is convex in $x$ and concave in $\lambda$ (see \cite{Xu:20}). It is noteworthy that when $\tau = 0$, the expression \eqref{eq:AugLag} simplifies to the classical augmented Lagrangian for problem \eqref{eq:prob}  as defined in \cite{Roc:76, Xu:20, JinWan:22}, i.e.,
\[ \mathcal{L}_{\rho, 0} (x;\lambda) := F(x) + \frac{1}{2\rho m}\sum_{j = 1}^m \left[ (\rho h_j(x) + \lambda^j)^2_+ - (\lambda^j)^2 \right].\]
\noindent Moreover,  we have the following expressions for the gradients of \eqref{eq:AugLag} w.r.t. $x$ and $\lambda$, for all $ j = 1:m$:
\begin{align} 
    & \nabla \psi^j_{\rho,\tau} (h_j(x); \lambda^j) = (\rho h_j(x) + (1-\tau) \lambda^j)_+\cdot \nabla h_j(x), \label{eq:nabla_x}\\
    & \nabla_\lambda \psi^j_{\rho,\tau} (h_j(x); \lambda^j) = \frac{1\!-\! \tau}{\rho} \left[(\rho h_j(x) \!+\! (1\!-\!\tau) \lambda^j)_+ \!-\! (1 \!-\!\tau)\lambda^j \right]\label{eq:nabla_lambda}
    = {(1\!-\! \tau)} \max \left( -\frac{(1\!-\!\tau)\lambda^j}{\rho} , h_j(x) \right).
\end{align}
\noindent The main advantage of having this perturbed augmented Lagrangian ($\tau > 0$) over standard augmented Lagrangian ($\tau = 0$) is that it allows to have dual updates well-behaved and easier to analyze (see next sections). 
Given $\rho > 0$ and $\tau \in [0,1)$, the dual function corresponding to \eqref{eq:AugLag} is defined as:
\begin{align}
    d_{\rho,\tau}(\lambda) = \min_{x \in \mathcal{Y}} \mathcal{L}_{\rho,\tau} (x; \lambda),
\end{align}
and any maximizer,  $\lambda^*$, of the dual problem $\max_{\lambda \geq 0} d_{\rho,\tau}(\lambda)$ is a dual   solution of \eqref{eq:prob}.  Now, we provide some basic properties for $\mathcal{L}_{\rho,\tau}$. First, we show that the gradient of the function $\psi^j_{\rho,\tau} (h_j(x); \lambda^j)$ w.r.t. $\lambda$ is bounded for all $j =1:m$.

\begin{lemma}\label{lemma_bdd_gradDual} Let us assume that the bounds from \eqref{eq:Bddness} hold. Then, choosing $\rho > 0, \tau \in [0,1)$ and $x \in \mathcal{Y}$, we have the following relation true:
    $$\frac{|\nabla_\lambda \psi^j_{\rho,\tau} (h_j(x); \lambda^j)|}{1-\tau}  = \left| \max \left( -\frac{(1- \tau)\lambda^j}{\rho}, h_{j}(x) \right) \right| \le M_h \;\;\;\;   \forall \lambda^j \ge 0, \;\;  j = 1:m.$$
\end{lemma}

\proof{Proof:} 
We prove this lemma by dividing it into two cases as follows:\\
    \textit{Case 1}: If $h_{j}(x) \ge - \frac{ (1- \tau)\lambda^{j}}{\rho}$, it follows that:
    \vspace{-0.3cm}
    $$\left| \max \left( -\frac{(1- \tau) \lambda^j}{\rho}, h_{j}(x) \right) \right| = |h_{j}(x)| \overset{\eqref{eq:Bddness}}{\le} M_h.$$
    \textit{Case 2}: If $h_{j}(x) < - \frac{ (1- \tau)\lambda^{j}}{\rho}$,  since $\lambda^j \ge 0$ for any $j =1:m$,   it follows that $h_{j}(x) < - \frac{(1- \tau) \lambda^{j}}{\rho} \le 0$, which gives us:
    $$\left| \max \left( -\frac{(1- \tau)\lambda^j}{\rho}, h_{j}(x) \right) \right| = \left| - \frac{(1- \tau)\lambda^j}{\rho} \right| \le |h_{j}(x)| \overset{\eqref{eq:Bddness}}{\le} M_h.$$
    Hence, combining both cases, we have the required result.  \Halmos
\endproof 

\medskip 

\noindent The following lemma shows that the function $\Psi_{\rho,\tau} (h(x);\lambda)$ is nonpositive. The proof of this lemma is similar to  Lemma 3.1 with $\tau = 0$ in \cite{Xu:20}, but here  adapted for \eqref{eq:AugLag} and it plays a key role in proving the convergence guarantees.

\begin{lemma}\label{lem:negative}
    Let $\rho > 0, \tau \in [0,1)$, then for any $x \in \mathcal{Y}$ such that $h_j(x) \le 0$, for all $j = 1:m$, and for any $\lambda \ge 0,$ it holds that:
    $$\Psi_{\rho,\tau} (h(x);\lambda) = \frac{1}{2\rho m}\sum_{j = 1}^m \left[(\rho h_j(x) + (1-\tau) \lambda^j)^2_+ - ((1-\tau) \lambda^j)^2 \right] \le 0.$$
\end{lemma}

\proof{Proof:} 
    We consider the following two cases.  \\ 
    \noindent \textit{Case 1}: If $\rho h_j(x) + (1-\tau)\lambda^j \le 0$, then since $\lambda \ge 0, \rho > 0$ and $\tau \in [0,1)$, we have:
    \[ \frac{1}{2\rho m}\sum_{j = 1}^m \left[(\rho h_j(x) + (1-\tau) \lambda^j)^2_+ - ((1-\tau) \lambda^j)^2 \right] = - \frac{1}{2\rho m}\sum_{j = 1}^m ((1 - \tau)\lambda^j)^2 \le 0. \]
    \textit{Case 2}: If $\rho h_j(x) + (1- \tau)\lambda^j > 0$, we have:
    \begin{align*}
        & \frac{1}{2\rho m}\sum_{j = 1}^m \left[(\rho h_j(x) + (1-\tau) \lambda^j)^2_+ - ((1-\tau)\lambda^j)^2 \right] = \frac{1}{2 \rho m}\sum_{j = 1}^m \left[\rho h_j(x) (\rho h_j(x) + 2 (1-\tau) \lambda^j) \right] \le 0,
    \end{align*}
    where the last inequality follows from   $x \in \mathcal{Y}$ and $h_j(x) \le 0$, for all $j = 1:m$. This proves the statement. \Halmos
\endproof 

\medskip 

\noindent Next lemma shows  the smoothness property of $\mathcal{L}_{\rho, \tau} (x; \lambda)$ w.r.t. $x$ for fixed $\lambda$.

\begin{lemma}\label{lemma_smooth}
Let Assumption \ref{Ass:Smooth} hold and the bounds from \eqref{eq:Bddness} are valid. Additionally,  let $\rho > 0, \tau \in [0,1), \lambda \in \mathbb{R}_+^m$ and $L_{\rho,\tau}^\lambda  = L_f + \rho B_h^2 + \left(\rho M_h + \frac{(1-\tau)}{m} \|\lambda\|_1 \right)L_h $. Then, the following is true:
    \[ \|\nabla \mathcal{L}_{\rho,\tau} (x; \lambda) -  \nabla \mathcal{L}_{\rho,\tau} (y; \lambda)\| \le L_{\rho, \tau}^\lambda \| x - y\| \quad  \forall x,y \in \mathcal{Y}. \]
\end{lemma}

\proof{Proof:} 
For any fixed $\lambda \in \mathbb{R}_+^m$, we have the following:
\begin{align*}
    \|\nabla \mathcal{L}_{\rho,\tau} (x; \lambda) -  \nabla \mathcal{L}_{\rho,\tau} (y; \lambda)\| & \!=\! \left\| \nabla F(x) + \frac{1}{m}\sum_{j = 1}^m \nabla \psi^j_{\rho,\tau} (h_j(x); \lambda^j) - \nabla F(y) - \frac{1}{m}\sum_{j = 1}^m \nabla \psi^j_{\rho,\tau} (h_j(y); \lambda^j) \right\|\\
    & \overset{\eqref{eq:nabla_x}}{\le} \|\nabla F(x) - \nabla F(y)\| + \frac{1}{m}\sum_{j = 1}^m (\rho h_{j}(x) + (1-\tau) 
    \lambda^{j})_+\|\nabla h_{j}(x) - \nabla h_{j}(y)\| \\
    & \quad + \frac{1}{m}\sum_{j = 1}^m |(\rho h_{j}(x) + (1- \tau)\lambda^{j})_+ - (\rho h_{j}(y) + (1- \tau)\lambda^{j})_+| \cdot\|\nabla h_{j}(y)\| \\
    & \le \|\nabla F(x) - \nabla F(y)\| + \frac{1}{m}\sum_{j = 1}^m |\rho h_{j}(x) + (1-\tau)\lambda^{j}| \|\nabla h_{j}(x) - \nabla h_{j}(y)\| \\
    & \quad + \rho \frac{1}{m}\sum_{j = 1}^m |h_{j}(x)  - h_{j}(y)| \cdot\|\nabla h_{j}(y)\|\\
    & \overset{\eqref{eq:Bddness}}{\le} \left(L_f + \rho B_h^2 + \left(\rho M_h + \frac{(1-\tau)}{m}\sum_{j = 1}^m|\lambda^{j}|\right)L_h \right) \|x - y\| = L_{\rho, \tau}^\lambda \|x - y\|,
\end{align*}
where the first inequality follows from the triangle inequality of the norm, the second inequality uses the fact that $(a)_+ \le |a|$ and the Lipschitz continuity of $(\cdot)_+$ function, in the third inequality we use the Lipschitz continuous gradient condition on $F$ and $h_{j}, \forall j = 1:m$, i.e., Assumption \ref{Ass:Smooth}. This concludes the proof.  \Halmos
\endproof 

\medskip 

\noindent This result implies the following descent lemma, see \cite{Nes:18}:
\begin{align}\label{eq:smoothPAL}
    |\mathcal{L}_{\rho,\tau} (y; \lambda_{k}) - \mathcal{L}_{\rho,\tau} (x; \lambda_{k}) - \langle \nabla \mathcal{L}_{\rho,\tau} (x; \lambda_{k}), y - x \rangle| \le \frac{L_{\rho, \tau}^\lambda}{2} \|x-y \|^2 \quad \forall x,y \in \mathcal{Y}.
\end{align}


\noindent One should notice here that the constant $L_{\rho, \tau}^\lambda$ depends on the dual multiplier $\lambda$, thus to ensure the well-definitness of this constant it is important to either assume or prove it to be bounded. In the next sections, we investigate this issue in detail. But before this, let us show that the strong duality holds for \eqref{eq:AugLag} under Assumption \ref{Ass:KKT}. Indeed, since we assume that the KKT conditions hold at a primal-dual solution $(x^*, \lambda^*)$ for \eqref{eq:prob}, then for any $\rho > 0$, $\tau \in [0,1)$, we have for $(\lambda^*)^j > 0$: 
$$(\rho h_j(x^*) + (1 - \tau)(\lambda^*)^j)_+ = \frac{(\lambda^*)^j (\rho h_j(x^*) + (1 - \tau)(\lambda^*)^j)_+}{(\lambda^*)^j} \overset{\eqref{eq:KKT3}}{=} (1 - \tau) (\lambda^*)^j. $$ 
A similar relation holds when $(\lambda^*)^j = 0$. Further, using \eqref{eq:nabla_x}, relation \eqref{eq:KKT1} becomes $0 \in \nabla \mathcal{L}_{\rho,\tau} (x^*;\lambda^*) + \mathcal{N}_\mathcal{Y} (x^*)$. Since $F$ and $h_j$'s are convex functions and $\mathcal{Y}$ is a convex set, it follows that  $x^*$ is a solution of $\min_{x \in \mathcal{Y}} \mathcal{L}_{\rho,\tau} (x;\lambda^*)$, which means $d_{\rho,\tau}(\lambda^*) = \mathcal{L}_{\rho,\tau} (x^*;\lambda^*)$. Further, using  \eqref{eq:KKT3} in \eqref{eq:AugLag}, we obtain:
\begin{align*}
    & \mathcal{L}_{\rho,\tau} (x^*;\lambda^*)  = F(x^*) + \frac{1-\tau}{2 \rho m}\sum_{j = 1}^m \left[(\lambda^*)^j(\rho h_j(x^*) + (1- \tau) (\lambda^*)^j)_+ - (1 -\tau)((\lambda^*)^j)^2 \right]  \overset{\eqref{eq:KKT3}}{=} F(x^*), 
\end{align*}
where we use the fact that $(a)_+^2 = a(a)_+$. Therefore, the strong duality holds, i.e., we have:
\[ d_{\rho,\tau}(\lambda^*) =  F(x^*). \]  
Further, we also have the following result that will be used later in our convergence analysis.

\begin{lemma}
    Let  $F$ and $h_j$'s, for all $ j=1:m$, be convex functions and also the set $\mathcal{Y}$ be convex. If a primal-dual solution $(x^*, \lambda^*)$ satisfies Assumption \ref{Ass:KKT}, then we have the following  true:
    \begin{align}\label{eq:KKTconv}
        F(x) - F(x^*) + \frac{1-\tau}{m} \sum_{j=1}^m (\lambda^*)^j h_j(x) \ge 0 \;\;\; \forall x \in \mathcal{Y}.
    \end{align}
\end{lemma}

\proof{Proof:} 
    Since we know that $(\lambda^*)^j \ge 0$ and $h_j$'s, for all $ j = 1:m$, are convex functions, we obtain:
    \begin{align}\label{eq:1SGDPA}
        (\lambda^*)^j (h_j (x) - h_j (x^*)) \ge \langle (\lambda^*)^j\nabla h_j (x^*), x - x^*\rangle.
    \end{align} 
    From \eqref{eq:KKT1}, we have:
    \[ -\frac{1-\tau}{m} \sum_{j = 1}^m (\lambda^*)^j \nabla h_j(x^*) \in \nabla F(x^*) + \mathcal{N}_\mathcal{Y} (x^*). \]
    From the convexity of $F$ and $\mathcal{Y}$,  we get for any $x \in \mathcal{Y}$:
    \begin{align*}
        F(x) & \ge  F(x^*) - \left\langle \frac{1-\tau}{m} \sum_{j = 1}^m (\lambda^*)^j \nabla h_j(x^*), x - x^* \right\rangle \overset{\eqref{eq:1SGDPA}}{\ge} F(x^*) - \frac{1-\tau}{m} \sum_{j = 1}^m (\lambda^*)^j (h_j (x) - h_j (x^*)) \\
        & \overset{\eqref{eq:KKT3}}{=} F(x^*) - \frac{1-\tau}{m} \sum_{j = 1}^m (\lambda^*)^j h_j (x),
    \end{align*}
where in the first inequality we used the definition of the indicator function and of the normal cone of $\mathcal{Y}$. Thus, we obtain the claimed result.  \Halmos
\endproof 


\section{Stochastic gradient descent and perturbed ascent method}
In this section, we present the Stochastic Gradient Descent and Perturbed Ascent (SGDPA) algorithm, which is a new stochastic primal-dual algorithm to solve problems of the form \eqref{eq:prob}. At each iteration, SGDPA takes a stochastic gradient of the perturbed augmented Lagrangian  w.r.t a single randomly chosen constraint to update the primal variable. Moreover, it also updates the dual variable in a randomly coordinate gradient ascent fashion  using the perturbation technique as in \cite{Lu:22, KosSha:11}. Hence, our new algorithm consists of the following iterative processes:
\begin{algorithm}[H]
\caption{\textbf{S}tochastic \textbf{G}radient \textbf{D}escent and \textbf{P}erturbed \textbf{A}scent (SGDPA) }\label{algorithmSGD}
\begin{algorithmic}[1]
\State  $\text{Choose} \; x_0 \in \mathcal{Y}, \lambda_0 \geq 0, \tau \in [0, 1), \; \text{stepsizes} \; \alpha_k>0$ and $\rho > 0$.
\State For {$k \geq 0$  update:}
\vspace{-0.3cm}
\begin{align}
    & \text{Sample $j_k \in [1:m]$ uniformly  at random and update primal variable:} \nonumber  \\
    & x_{k+1} = \Pi_{\mathcal{Y}} \left(x_{k} - \alpha_k \left( \nabla F(x_k) + \nabla \psi^{j_k}_{\rho,\tau} (h_{j_k}(x_k); \lambda^{j_k}_{k}) \right) \right)\label{eq:primal_update}\\
    & \text{Sample $\bar{j}_k \in [1:m]$ uniformly  at random ({independent on $j_k$}) and then update dual variable:} \nonumber\\
    & \lambda_{k+1}^j = \begin{cases}
    \lambda_k^j, & \text{if $j \neq \bar{j}_k$}\\
    (1-\tau)\lambda_k^j + \rho \max \left( -\frac{(1 - \tau)\lambda_k^j}{\rho}, h_j(x_{k+1}) \right), & \text{if $j = \bar{j}_k$}. \end{cases}\label{eq:dual_update}
\end{align}
\end{algorithmic}
\end{algorithm}

\noindent   The computational simplicity of the SGDPA subproblem, which uses only gradient information and considers only one constraint at a time, makes our algorithm suitable for problems with many functional constraints.  Indeed,  the primal update \eqref{eq:primal_update} can be expressed equivalently~as:
\[  x_{k+1} = \argmin_{x\in \mathcal{Y}} \langle \nabla F(x_k) + \nabla \psi^{j_k}_{\rho,\tau} (h_{j_k}(x_k); \lambda^{j_k}_{k}), x - x_k\rangle + \frac{1}{2\alpha_k} \|x - x_k\|^2. \]
\noindent  Hence, the primal update rule is nothing but the minimization of the linearized function $F(x) +  \psi^{j_k}_{\rho,\tau} (h_{j_k}(x); \lambda^{j_k}_{k})$ at the current iterate  with a proximal term added to it. The dual update \eqref{eq:dual_update} is equivalent to a random coordinate gradient ascent step:
\begin{align}\label{eq:fullupdate}
 \lambda_{k+1}  &= \begin{cases}
        \lambda_k^j, & \text{if $j \neq \bar{j}_k$}\\
    \argmax_{\lambda \in \mathbb{R}_+} \left\langle \frac{1}{1 - \tau}\nabla_\lambda \mathcal{L}^{j}_{\rho,\tau} (x_{k+1};\lambda_k), \lambda - \lambda^{j}_k \right\rangle - \frac{1-\tau}{2 \rho m} \|\lambda-\lambda^{j}_k\|^2- \frac{\tau}{2\rho m}\|\lambda\|^2, & \text{if $j = \bar{j}_k$}
    \end{cases} \nonumber \\
    & = \lambda_k + e_{\bar{j}_k} \odot \left( \frac{\rho m}{(1 - \tau)} \nabla_{\lambda} \mathcal{L}_{\rho,\tau} (x_{k+1}; \lambda_k) - \tau \lambda_k \right),
\end{align}
where $\odot$ denotes the componentwise product between two vectors, the index $\bar{j}_k$ is chosen uniformly at random from $[1:m]$ and independent from $j_k$ and $e_{\bar{j}_k}\in \mathbb{R}^m$ denotes the vector with all entries zeros except the $\bar{j}_k^{th}$ entry, which is one. The perturbation parameter $\tau$ introduces negative curvature to the maximization problem, ensuring that the dual update is well-behaved. Also, as one will notice later in Section \ref{ch04:S2sss} (see Lemma \ref{lem:dualbound}) using this  dual perturbation technique we can prove  boundedness of the dual iterates. Note that  our algorithm is similar to the one in \cite{JinWan:22} when $\tau = 0$  and  the minibatch size for the constraints  is set to $1$ or  to  PDSG algorithm in \cite{Xu:20} with $\tau = 0$. However, SGDPA algorithm  differentiates itself from \cite{JinWan:22, Xu:20} by  using two distinct uniformly random samples for the primal and dual updates from the set of indices $[1:m]$; SGDPA is also based on a general perturbed augmented Lagrangian framework as $\tau \in [0, 1)$;   SGDPA uses  $x_{k+1}$ instead of $x_k$  used in PDSG from  \cite{Xu:20} to update the dual multipliers; the algorithm  in \cite{JinWan:22} needs to work with the full set of constraints at each iteration in order to get optimal performance,  while SGDPA uses only one constraint at a time. It is important to notice that our strategy of separate samplings for the primal and dual updates allows us to update the dual variable at the current state $x_{k+1}$, as it will be evident later in the proof of Lemma \ref{lemma_dualupdate}. 

\medskip 

\noindent Further, let us define the filtration:
\[ \mathcal{F}_{[k]} := \{ j_0,..., j_k, \bar{j}_0,..., \bar{j}_k \}.  \]
Also, one can easily notice that since $j_k$ is chosen uniformly at random,   $\nabla \psi^{j_k}_{\rho,\tau} (h_{j_k}(x_k); \lambda^{j_k}_k) $ is an unbiased estimator of $\nabla \Psi_{\rho,\tau} (h(x_k), \lambda_k)$, i.e.: 
$$\mathbb{E}_{j_k} [\nabla \psi^{j_k}_{\rho,\tau} (h_{j_k}(x_k); \lambda^{j_k}_k) | \mathcal{F}_{[k-1]}] = \nabla \Psi_{\rho,\tau} (h(x_k), \lambda_k).$$
Similarly, $\psi^{j_k}_{\rho,\tau} (h_{j_k}(x_k); \lambda^{j_k}_k) $ is an unbiased estimator of $\Psi_{\rho,\tau} (h(x_k), \lambda_k)$, i.e.: 
$$\mathbb{E}_{j_k} [\psi^{j_k}_{\rho,\tau} (h_{j_k}(x_k); \lambda^{j_k}_k) | \mathcal{F}_{[k-1]}] = \Psi_{\rho,\tau} (h(x_k), \lambda_k).$$


\section{Convergence analysis of SGDPA}
\noindent This section presents the convergence results for the SGDPA algorithm under two different strategies for the dual update: (i) when the dual iterates, $\lambda_k$, are assumed to be bounded; and (ii) when we use the dual perturbation update technique to prove that the iterates $\lambda_k$ are bounded. In both scenarios, we achieve  sublinear rates of convergence when the objective is  convex or strongly convex. Before deriving the main results of this section, we provide some preliminary results.  First, we show that the dual variables generated by \eqref{eq:dual_update} are non-negative. 
\begin{lemma}\label{lem:nonnegDual}
Let $\rho >0 $ and $ \tau \in [0,1)$. Then,  the dual iterates generated by SGDPA satisfy $\lambda_k \geq 0$ for all $k \ge 0$. 
\end{lemma}

\proof{Proof:} 
We use induction to prove the result. First, the result holds when $k = 0$ due to our choice $\lambda_0 \in \mathbb{R}^m_+$. Assume the statement holds for $k > 0$, i.e., $\lambda_k \in \mathbb{R}^m_+$. Then, at iteration $k$, it follows from the dual update \eqref{eq:dual_update} in  SGDPA  with given choices of $\tau, \rho$, for any $\bar{j}_k$ sampled uniformly at random from $[1:m]$, we have:
    \begin{align*}
    \lambda_{k+1}^{j}
    \begin{cases}
        = \lambda_{k}^{j} \in \mathbb{R}_+ & \text{ when } j \neq \bar{j}_k,\\
        \ge (1 - \tau)\lambda_k^{j} - \rho \frac{(1 - \tau)\lambda_k^{j}}{\rho} = 0  & \text{ when } j = \bar{j}_k.
    \end{cases}
    \end{align*} 
   Thus $\lambda_{k+1} \in \mathbb{R}^m_+$. This completes the proof. \Halmos
\endproof 

\medskip 

\noindent In the next lemma, we provide the descent type nature of the perturbed augmented Lagrangian function, see  \eqref{eq:AugLag},  along the primal iterates.
\begin{lemma}\label{lem:x_kx_k+1}
Let Assumption \ref{Ass:Smooth} hold and, additionally,  the bounds from \eqref{eq:Bddness} be true. Further, let $\rho > 0, \tau \in [0, 1)$, and choose a positive stepsize sequence $\{\alpha_k\}_{k\ge 0}$. Then, the following relation holds in expectation   for the sequences generated by SGDPA:
    \begin{align*}
         & \mathbb{E}[\mathcal{L}_{\rho,\tau} (x_{k+1}; \lambda_k)] \\
         & \le \mathbb{E}[\mathcal{L}_{\rho,\tau} (x_k; \lambda_k)] - \mathbb{E}\left[ \frac{2 - \alpha_k L_{\rho,\tau}^{\lambda_k}}{2\alpha_k} \|x_{k+1} - x_k\|^2\right] + \mathbb{E}\left[ 2 \alpha_k \left( \rho^2 M_h^2 + (1-\tau)^2 \|\lambda_k\|^2 \right)B_h^2 \right] \quad \forall k \geq 0.
    \end{align*}
\end{lemma}

\proof{Proof:} 
From the optimality condition of the update \eqref{eq:primal_update}, we have:
\begin{align*}
    & \left\langle x_{k+1} - x_k + \alpha_k  \left( \nabla F(x_k) + \nabla \psi^{j_k}_{\rho,\tau} (h_{j_k}(x_k); \lambda^{j_k}_{k})\right), x_{k+1} - x_k \right\rangle \le 0,
    \end{align*}
   which implies  
    \begin{align}\label{eq:x_kx_k+1} 
    \frac{1}{\alpha_k}\|x_{k+1} - x_k\|^2 + \langle \nabla F(x_k) + \nabla \psi^{j_k}_{\rho,\tau} (h_{j_k}(x_k); \lambda^{j_k}_{k}), x_{k+1} - x_k \rangle\le 0.
\end{align}
From \eqref{eq:smoothPAL}, we get:
\begin{align}\label{eq:11SGDPA}
    &  \mathcal{L}_{\rho,\tau} (x_{k+1}; \lambda_k) \le \mathcal{L}_{\rho,\tau} (x_k; \lambda_k) + \langle \nabla \mathcal{L}_{\rho,\tau} (x_k; \lambda_k), x_{k+1} - x_k\rangle \!+\!\frac{L_{\rho,\tau}^{\lambda_k}}{2} \|x_{k+1} - x_k\|^2 \\
    & \overset{\eqref{eq:x_kx_k+1}}{\le} \mathcal{L}_{\rho,\tau} (x_k; \lambda_k) + \langle \nabla \Psi_{\rho,\tau} (h(x_k); \lambda_k) - \nabla \psi^{j_k}_{\rho,\tau} (h_{j_k}(x_k); \lambda^{j_k}_{k}), x_{k+1} - x_k\rangle  - \frac{2 - \alpha_k L_{\rho,\tau}^{\lambda_k}}{2\alpha_k} \|x_{k+1} - x_k\|^2. \nonumber
\end{align}
Denote, $\Tilde{x}_{k+1} = \Pi_\mathcal{Y} \left(x_k - \alpha_k \left( {\nabla F (x_k)} + \nabla \Psi_{\rho,\tau} (h(x_k); \lambda_k)\right) \right).$ By the nonexpansiveness of the projection operator, we have:
\begin{align*}
    \|\Tilde{x}_{k+1} - x_{k+1}\| \le \alpha_k\| \nabla \Psi_{\rho,\tau} (h(x_k); \lambda_k) - \nabla \psi^{j_k}_{\rho,\tau} (h_{j_k}(x_k); \lambda^{j_k}_{k})\|.
\end{align*}
Now, using it in \eqref{eq:11SGDPA}, we get:
\begin{align*}
    \mathcal{L}_{\rho,\tau} (x_{k+1}; \lambda_k)  & \le \mathcal{L}_{\rho,\tau} (x_k; \lambda_k) + \langle \nabla \Psi_{\rho,\tau} (h(x_k); \lambda_k) - \nabla \psi^{j_k}_{\rho,\tau} (h_{j_k}(x_k); \lambda^{j_k}_{k}), x_{k+1} - \Tilde{x}_{k+1}\rangle\\
    &\quad + \langle \nabla \Psi_{\rho,\tau} (h(x_k); \lambda_k) - \nabla \psi^{j_k}_{\rho,\tau} (h_{j_k}(x_k); \lambda^{j_k}_{k}), \Tilde{x}_{k+1} - x_k\rangle  - \frac{2 - \alpha_k L_{\rho,\tau}^{\lambda_k}}{2\alpha_k} \|x_{k+1} - x_k\|^2\\
    & \le \mathcal{L}_{\rho,\tau} (x_k; \lambda_k) + \alpha_k \| \nabla \Psi_{\rho,\tau} (h(x_k); \lambda_k) - \nabla \psi^{j_k}_{\rho,\tau} (h_{j_k}(x_k); \lambda^{j_k}_{k})\|^2\\
    & \quad + \langle \nabla \Psi_{\rho,\tau} (h(x_k); \lambda_k) - \nabla \psi^{j_k}_{\rho,\tau} (h_{j_k}(x_k); \lambda^{j_k}_{k}), \Tilde{x}_{k+1} - x_k\rangle - \frac{2 - \alpha_k L_{\rho,\tau}^{\lambda_k}}{2\alpha_k} \|x_{k+1} - x_k\|^2,
\end{align*}
where in the last inequality we used the Cauchy-Schwartz inequality. After taking expectation conditioned on $j_k$, we obtain {(since the third term will be 0)}:
\begin{align*}
     \mathbb{E}_{j_k}[\mathcal{L}_{\rho,\tau} (x_{k+1}; \lambda_k)|\mathcal{F}_{[k-1]}]
    & \le \mathcal{L}_{\rho,\tau} (x_k; \lambda_k) - \frac{2 - \alpha_k L_{\rho,\tau}^{\lambda_k}}{2\alpha_k} \mathbb{E}_{j_k}[\|x_{k+1} - x_k\|^2 |\mathcal{F}_{[k-1]}]\\
    & \quad + \alpha_k  \mathbb{E}_{j_k}[ \| \nabla \Psi_{\rho,\tau} (h(x_k); \lambda_k) - \nabla \psi^{j_k}_{\rho,\tau} (h_{j_k}(x_k); \lambda^{j_k}_{k})\|^2 |\mathcal{F}_{[k-1]}]\\
    & \le \mathcal{L}_{\rho,\tau} (x_k; \lambda_k) - \frac{2 - \alpha_k L_{\rho,\tau}^{\lambda_k}}{2\alpha_k} \mathbb{E}_{j_k}[\|x_{k+1} - x_k\|^2 |\mathcal{F}_{[k-1]}]\\
    & \quad + \alpha_k  \mathbb{E}_{j_k}[ \| \nabla \psi^{j_k}_{\rho, \tau} (h_{j_k}(x_k); \lambda^{j_k}_{k})\|^2 |\mathcal{F}_{[k-1]}]  \\
    & \overset{\eqref{eq:nabla_x}}{=} \mathcal{L}_{\rho,\tau} (x_k; \lambda_k) - \frac{2 - \alpha_k L_{\rho,\tau}^{\lambda_k}}{2\alpha_k} \mathbb{E}_{j_k}[\|x_{k+1} - x_k\|^2 |\mathcal{F}_{[k-1]}]\\
    & \quad + \alpha_k  \mathbb{E}_{j_k}[ \|(\rho h_{j_k} (x_k) + (1-\tau) \lambda_k^{j_k})_+ \nabla h_{j_k} (x_k) \|^2 |\mathcal{F}_{[k-1]}]\\
    & \overset{\eqref{eq:Bddness}}{\le} \mathcal{L}_{\rho,\tau} (x_k; \lambda_k) - \frac{2 - \alpha_k L_{\rho,\tau}^{\lambda_k}}{2\alpha_k} \mathbb{E}_{j_k}[\|x_{k+1} - x_k\|^2 |\mathcal{F}_{[k-1]}]\\
    & \quad + 2 \alpha_k \left( \rho^2 M_h^2 + (1-\tau)^2 \|\lambda_k\|^2 \right) B_h^2,
\end{align*}
where the second inequality follows from the fact that $\mathbb{E}[\|X - \mathbb{E}[X]\|^2] = \mathbb{E}[\|X - Y\|^2] - \mathbb{E}[\|Y - \mathbb{E}[X]\|^2]$, for any random variable $X$ and for any $Y$ and in the last inequality we use the inequality $(a+b)_+^2 \le |a+b|^2 \le 2|a|^2 + 2|b|^2$ together with the fact that the index $j_k$ is chosen uniformly at random. Now after taking the full expectation we obtain the claimed result. \Halmos
\endproof 

\medskip 

\noindent The following lemma gives a descent type relation between $\lambda_{k+1}$ and $\lambda_k$ using the dual update \eqref{eq:dual_update}. Before proving this relation let us denote: $$\Delta_{\rho,\tau} (x_{k+1}; \lambda_k) = \frac{1}{(1 - \tau)}\nabla_\lambda \mathcal{L}_{\rho,\tau} (x_{k+1} ; \lambda_k) - \frac{\tau}{\rho m} \lambda_k.$$

\begin{lemma}\label{lemma_dualupdate}
For any deterministic or stochastic  $ \lambda \ge 0$ (independent from $\bar{j}_k$ at $k^{\text{th}}$ iteration),  $\rho > 0, \tau \in [0, 1)$, it holds the following for the iterates generated by  SGDPA algorithm:
\begin{align*}
    &\mathbb{E}\left[\frac{1 - \tau}{m} \sum_{j=1}^m \lambda^j h_j (x_{k+1}) - \Psi_{\rho, \tau} (h(x_{k+1}); \lambda_k ) \right] \\
    & \le \frac{1}{2\rho m} \left(\mathbb{E}[\|\lambda_k - \lambda\|^2] - \mathbb{E}[\|\lambda_{k+1} - \lambda\|^2]\right) + \frac{\tau}{\rho m}\mathbb{E} [\|\lambda\|^2]   \quad \forall k \geq 0.
\end{align*}
\end{lemma}

\proof{Proof:} 
Let us denote $J_k^+ = \{ j\in [m] : \rho h_j(x_{k+1}) + (1 - \tau) \lambda_k^j \ge 0 \}$  and $J_k^- = [m]\backslash J_k^+.$
Then, for any $ \lambda \ge 0$, we have:
\begin{align*}
    & \frac{1 - \tau}{m} \sum_{j=1}^m \lambda^j h_j (x_{k+1}) - \Psi_{\rho, \tau} (h(x_{k+1}); \lambda_k )\\
    & = \frac{1 - \tau}{m} \sum_{j=1}^m \lambda^j h_j (x_{k+1}) - \frac{1}{2\rho m} \sum_{j =1}^m [(\rho h_j(x_{k+1}) + (1 - \tau) \lambda^{j}_{k})^2_+ - ((1 - \tau) \lambda^{j}_{k})^2] \\
    & = - \frac{1}{2 \rho m} \left[\sum_{j\in J_k^+} (\rho h_j(x_{k+1}))^2 + 2 \rho h_j(x_{k+1}) (1 - \tau) \lambda_k^j - 2 \rho h_j(x_{k+1}) (1 - \tau) \lambda^j + (\tau\lambda_k^j)^2 - (\tau\lambda_k^j)^2 \right] \\
    & \quad + \frac{1}{2 \rho m}\left[\sum_{j \in J_k^-} ((1 - \tau)\lambda_k^j)^2 + 2 \rho h_j(x_{k+1}) (1 - \tau)\lambda^j + (\lambda_k^j)^2 - (\lambda_k^j)^2 \right] \\
    & = -\frac{1}{2 \rho m} \left[\sum_{j\in J_k^+} (\rho h_j(x_{k+1}) - \tau \lambda_k^j)^2 + 2\rho h_j(x_{k+1})\lambda_k^j - (\tau\lambda_k^j)^2 - 2 \rho h_j(x_{k+1}) (1 - \tau)\lambda^j \right] \\
    & \quad + \frac{1}{2 \rho m}\left[\sum_{j \in J_k^-} ((1 - \tau)\lambda_k^j)^2 + 2 \rho h_j(x_{k+1}) (1 - \tau)\lambda^j + (\lambda_k^j)^2 - (\lambda_k^j)^2 \right] \\
    & = - \frac{1}{2\rho m} \left[\sum_{j\in J_k^+} (\rho h_j(x_{k+1}) - \tau \lambda_k^j)^2  + \sum_{j\in J_k^-} (\lambda_k^j)^2 \right] \\
    & \quad - \sum_{j\in J_k^+} \left[\frac{1}{m} h_j(x_{k+1})\lambda_k^j - \frac{(\tau\lambda_k^j)^2}{2\rho m} - \frac{(1 - \tau)}{m} h_j(x_{k+1}) \lambda^j \right] \\
    & \quad + \sum_{j \in J_k^-} \left[\frac{((1 - \tau)\lambda_k^j)^2}{\rho m} + \frac{(1 - \tau)}{m} h_j(x_{k+1})\lambda^j + \frac{1 - (1 - \tau)^2}{2\rho m}(\lambda_k^j)^2 \right] \\
    & = - \frac{1}{2\rho m} \left[\sum_{j\in J_k^+} (\rho h_j(x_{k+1}) - \tau \lambda_k^j)^2  + \sum_{j\in J_k^-} (\lambda_k^j)^2 \right] \\
    & \quad - \sum_{j\in J_k^+} \left[\left( \frac{1}{m} h_j(x_{k+1}) - \frac{\tau}{\rho m} \lambda_k^j \right)(\lambda_k^j - \lambda^j) + \frac{\tau(2 - \tau)}{2 \rho m}((\lambda_k^j)^2 - 2 \lambda_k^j\lambda^j) \right] \\
    & \quad + \sum_{j \in J_k^-} \left[\left(\frac{1}{\rho m} \lambda_k^j \right)(\lambda_k^j - \lambda^j) - \frac{\tau(2 - \tau)}{2 \rho m}((\lambda_k^j)^2 - 2 \lambda_k^j\lambda^j) \right] \\
    & \quad - \sum_{j\in J_k^+} \frac{\tau \lambda^j}{\rho m} (\rho h_j(x_{k+1}) + (1 - \tau)\lambda_k^j) + \sum_{j \in J_k^-}  \frac{(1 - \tau)\lambda^j}{\rho m} (\rho h_j(x_{k+1})+ (1 - \tau)\lambda_k^j). 
\end{align*}
\noindent Further, using the fact that $\rho h_j(x_{k+1}) + (1 - \tau) \lambda_k^j \ge 0$ for $j \in J_k^+$ and $\rho h_j(x_{k+1}) + (1 - \tau) \lambda_k^j < 0$ for $j \in J_k^-$, we get:
\begin{align*}
    & \frac{1 - \tau}{m} \sum_{j=1}^m \lambda^j h_j (x_{k+1}) - \Psi_{\rho, \tau} (h(x_{k+1}); \lambda_k )\\
    & \!\!\le \!- \frac{1}{2\rho m} \left[\sum_{j\in J_k^+} (\rho h_j(x_{k+1}) - \tau \lambda_k^j)^2  + \sum_{j\in J_k^-} (\lambda_k^j)^2 + \tau(2 - \tau) \sum_{j=1}^m (\lambda_k^j)^2 - 2 \lambda_k^j \lambda^j \right] \\
    & \quad - \left[\sum_{j\in J_k^+} \left( \frac{1}{m} h_j(x_{k+1}) - \frac{\tau}{\rho m} \lambda_k^j \right)(\lambda_k^j - \lambda^j) - \sum_{j \in J_k^-} \left(\frac{1}{\rho m} \lambda_k^j \right)(\lambda_k^j - \lambda^j) \right].
\end{align*}
Now, noticing that:
\begin{align*}
    & - \frac{1}{2\rho m} \left[\sum_{j\in J_k^+} (\rho h_j(x_{k+1}) - \tau \lambda_k^j)^2  + \sum_{j\in J_k^-} (\lambda_k^j)^2 \right] = - \frac{1}{2 \rho} \mathbb{E}_{\bar{j}_k} [\|\lambda_{k+1} - \lambda_k\|^2 | \mathcal{F}_{[k-1]}\cup j_k],
\end{align*}
the relation $- \|b\|^2 + 2 \langle a, b \rangle \le \|a\|^2$, the fact that $(2 -\tau) < 2$ and recalling that  $\Delta_{\rho,\tau} (x_{k+1}; \lambda_k) = \frac{1}{(1 - \tau)}\nabla_\lambda \mathcal{L}_{\rho,\tau} (x_{k+1} ; \lambda_k) - \frac{\tau}{\rho m} \lambda_k$, we get:
\begin{align}\label{eq:way}
    & \frac{1 - \tau}{m} \sum_{j=1}^m \lambda^j h_j (x_{k+1}) - \Psi_{\rho, \tau} (h(x_{k+1}); \lambda_k ) \\
    & \le - \frac{1}{2 \rho} \mathbb{E}_{\bar{j}_k} [\|\lambda_{k+1} - \lambda_k \|^2| \mathcal{F}_{[k-1]}\cup j_k] - \mathbb{E}_{\bar{j}_k} [\left\langle \Delta_{\rho,\tau} (x_{k+1}; \lambda_k), \lambda_k -\lambda \right\rangle | \mathcal{F}_{[k-1]}\cup j_k] + \frac{\tau}{\rho m} \|\lambda\|^2.\nonumber
\end{align}
Now, we bound the second term in the  right hand side of \eqref{eq:way}. From \eqref{eq:fullupdate}, we have:
\begin{align*}
   &  \left\langle \Delta_{\rho,\tau} (x_{k+1}; \lambda_k), \lambda_k -\lambda \right\rangle \\
    & = \left\langle e_{\bar{j}_k} \odot \left(\Delta_{\rho,\tau} (x_{k+1}; \lambda_k)\right), \lambda_k -\lambda \right\rangle - \left\langle e_{\bar{j}_k} \odot \left(\Delta_{\rho,\tau} (x_{k+1}; \lambda_k)\right)  - \left(\Delta_{\rho,\tau} (x_{k+1}; \lambda_k)\right), \lambda_k -\lambda \right\rangle\\
    & \overset{\eqref{eq:fullupdate}}{=} \frac{1}{\rho m} \langle \lambda_{k+1} - \lambda_k, \lambda_k - \lambda \rangle - \left\langle e_{\bar{j}_k} \odot \left(\Delta_{\rho,\tau} (x_{k+1}; \lambda_k)\right)  - \left(\Delta_{\rho,\tau} (x_{k+1}; \lambda_k)\right), \lambda_k -\lambda \right\rangle\\
    & = \!\frac{1}{2\rho m}[ \|\lambda_{k+1} \!-\!\lambda \|^2 \!-\! \|\lambda_k - \lambda\|^2 - \|\lambda_{k+1} - \lambda_k \|^2] - \left\langle e_{\bar{j}_k} \odot \Delta_{\rho,\tau} (x_{k+1}; \lambda_k)- \Delta_{\rho,\tau} (x_{k+1}; \lambda_k), \lambda_k -\lambda \right\rangle,
\end{align*}
where in the last equality we use the identity $2 \langle a - b, b - c \rangle = \|a - c \|^2 - \|b - c\|^2 - \|a - b\|^2$, for any $a,b,c \in \mathbb{R}^m$. Since  $\lambda \ge 0$ is assumed  deterministic or  stochastic (but independent on $\bar{j}_k$ at $k^{\text{th}}$ iteration), then taking the conditional expectation w.r.t. $\bar{j}_k$, we get:
\begin{align*}
   &  \mathbb{E}_{\bar{j}_k} [\left\langle \Delta_{\rho,\tau} (x_{k+1}; \lambda_k), \lambda_k -\lambda \right\rangle | \mathcal{F}_{[k-1]}\cup j_k] = \frac{1}{2\rho m} \mathbb{E}_{\bar{j}_k} [ \|\lambda_{k+1} -\lambda \|^2 - \|\lambda_k \!-\! \lambda\|^2 \!-\! \|\lambda_{k+1} - \lambda_k \|^2| \mathcal{F}_{[k-1]}\cup j_k].
\end{align*}
Using this in \eqref{eq:way}, we obtain:
\begin{align*}
    & \frac{1 - \tau}{m} \sum_{j=1}^m \lambda^j h_j (x_{k+1}) - \Psi_{\rho, \tau} (h(x_{k+1}); \lambda_k )  \\
    & {\le \frac{1}{2\rho m} \mathbb{E}_{\bar{j}_k} [ \|\lambda_k \!-\! \lambda\|^2 - \|\lambda_{k+1} \!-\!\lambda \|^2 | \mathcal{F}_{[k-1]}\cup j_k] - \frac{1}{2 \rho} \left( 1 - \frac{1}{m} \right) \mathbb{E}_{\bar{j}_k} [\|\lambda_{k+1} - \lambda_k \|^2| \mathcal{F}_{[k-1]}\cup j_k] + \frac{\tau}{\rho m} \|\lambda\|^2\nonumber}\\
    & {\overset{m \ge 1}{\le}\frac{1}{2\rho m} \mathbb{E}_{\bar{j}_k} [ \|\lambda_k - \lambda\|^2 - \|\lambda_{k+1} - \lambda \|^2 | \mathcal{F}_{[k-1]}\cup j_k] + \frac{\tau}{\rho m} \|\lambda\|^2}.
\end{align*}

\noindent After taking full expectation, we obtain the desired result.  \Halmos
\endproof 

\medskip 

\noindent The next theorem provides the main recurrence for the convergence analysis of SGDPA algorithm. Let us first introduce the following notation: $$B_{\tau, \lambda_k}^2 = \left(2 B^2_f + 8 (\rho^2 M_h^2 + (1 - \tau)^2 \|\lambda_k\|^2) B_h^2 \right).$$

\begin{theorem}\label{lemma_Main}
Let Assumptions \ref{Ass:KKT}, \ref{Ass:Smooth} and \ref{Ass:StrongConv} hold and, additionally,  the bounds from \eqref{eq:Bddness} be valid. Further,  let $ \rho > 0, \tau\in [0, 1)$ and choose a non-increasing stepsize $\alpha_k \in \left(0, \frac{2}{L_{\rho,\tau}^{\lambda_k}} \right]$, for all $k\ge 0$. Then, for any $x^* \in \mathcal{X}^*$ and $ \lambda \ge 0$, we have the following recurrence true for the iterates generated by  SGDPA algorithm:
\begin{align*}
    & \mathbb{E} [ \|x_{k+1} - x^* \|^2 ] + \frac{\alpha_k}{\rho m} \mathbb{E}[\|\lambda_{k+1} - \lambda\|^2]  + 2\alpha_k \mathbb{E} \left[F(x_{k+1}) - F^* +\frac{1 - \tau}{m} \sum_{j=1}^m \lambda^j h_j (x_{k+1})\right] \\
    & \le {(1 \!-\! \mu \alpha_k)} \mathbb{E}[\|x_k - x^*\|^2] + \frac{\alpha_k}{\rho m} \mathbb{E}[\|\lambda_k - \lambda\|^2]   + 2\alpha_k \frac{\tau}{\rho m}\mathbb{E}[\|\lambda\|^2]  +\mathbb{E}\left[ \alpha_k^2 B_{\tau, \lambda_k}^2 \right].
\end{align*}
\end{theorem}

\proof{Proof:} 
 Recall first, the following property of the projection onto a convex set $\mathcal{Y}$:
    \[ \|\Pi_\mathcal{Y} (v) - y\|^2 \le \|v-y\|^2 - \|\Pi_\mathcal{Y} (v) - v\|^2 \quad  \forall v \in \mathbb{R}^n \text{ and } y\in \mathcal{Y}. \]
Using this relation with $v = x_k - \alpha_k (\nabla F(x_k) + \nabla \psi^{j_k}_{\rho,\tau} (h_{j_k}(x_k); \lambda^{j_k}_{k})) \in \mathbb{R}^n$ and  $y = x^* \in \mathcal{Y}$, we get:
\begin{align*}
&\|x_{k+1} - x^*\|^2 \le \|x_k - x^* - \alpha_k (\nabla F(x_k) + \nabla \psi^{j_k}_{\rho,\tau} (h_{j_k}(x_k); \lambda^{j_k}_{k}))\|^2 \\
& = \|x_k - x^*\|^2 + \alpha_k^2 \|\nabla F(x_k) + \nabla \psi^{j_k}_{\rho,\tau} (h_{j_k}(x_k); \lambda^{j_k}_{k})\|^2  + 2\alpha_k \langle \nabla F(x_k) + \nabla \psi^{j_k}_{\rho,\tau} (h_{j_k}(x_k); \lambda^{j_k}_{k}), x^* - x_k \rangle\\
& \overset{\eqref{eq:nabla_x}, \text{ Ass.} \ref{Ass:StrongConv}}{\le} {(1 - \mu \alpha_k)}\|x_k - x^*\|^2 +\! 2 \alpha_k^2 (\|\nabla F(x_k)\|^2 + (\rho h_{j_k} (x_k) + (1-\tau) \lambda_k^{j_k})_+^2 \| \nabla h_{j_k} (x_k) \|^2 ) \\
& \qquad \qquad  + 2 \alpha_k (F(x^*) + \psi^{j_k}_{\rho,\tau} (h_{j_k}(x^*); \lambda^{j_k}_{k}) - F(x_k) - \psi^{j_k}_{\rho,\tau} (h_{j_k}(x_k); \lambda^{j_k}_{k})) \\
& \le  {(1 - \mu \alpha_k)} \|x_k - x^*\|^2 + 2 \alpha_k^2 (\|\nabla F(x_k)\|^2 + 2(|\rho h_{j_k} (x_k)|^2 + |(1-\tau) \lambda_k^{j_k}|^2) \|\nabla h_{j_k} (x_k) \|^2) \\
& \qquad + 2 \alpha_k (F(x^*) + \psi^{j_k}_{\rho,\tau} (h_{j_k}(x^*); \lambda^{j_k}_{k}) - F(x_k) - \psi^{j_k}_{\rho,\tau} (h_{j_k}(x_k); \lambda^{j_k}_{k})) \\
& \overset{\eqref{eq:Bddness}}{\le} {(1 - \mu \alpha_k)} \|x_k - x^*\|^2 + 2 \alpha_k^2 \left(B_F^2 + 2(\rho^2 M_h^2 + (1 - \tau)^2 |\lambda_k^{j_k}|^2) B_h^2\right)\\
& \qquad + 2 \alpha_k (F(x^*) + \psi^{j_k}_{\rho,\tau} (h_{j_k}(x^*); \lambda^{j_k}_{k}) - F(x_k) - \psi^{j_k}_{\rho,\tau} (h_{j_k}(x_k); \lambda^{j_k}_{k})),
\end{align*}
where the second inequality follows from the strong convexity of $F$ (see Assumption \ref{Ass:StrongConv}), the convexity of $h_j$'s and the relation $\|a+b\|^2 \le 2 \|a\|^2 + 2\|b\|^2$, for any $a,b \in \mathbb{R}^n$. Now, after taking conditional expectation w.r.t. $j_k$, we get:
\begin{align*}
    & \mathbb{E}_{j_k} [ \|x_{k+1} - x^* \|^2| \mathcal{F}_{[k-1]} ] \\
    & \le {(1 - \mu \alpha_k)}\|x_k - x^*\|^2 \!-\! 2\alpha_k (\mathcal{L}_{\rho,\tau} (x_k; \lambda_{k}) \!-\! \mathcal{L}_{\rho,\tau} (x^*; \lambda_{k}))  + 2\alpha_k^2 \left(B^2_f \!+\! 2 (\rho^2 M_h^2 \!+\! (1 \!-\! \tau)^2 \|\lambda_k\|^2) B_h^2\right).
\end{align*}
After taking full expectation, using the result from Lemma \ref{lem:x_kx_k+1},  noticing that $\alpha_k \le 2/L_{\rho,\tau}^{\lambda_k}$ and  $B_{\tau, \lambda_k}^2 = \left(2 B^2_f + 8 (\rho^2 M_h^2 \!+\! (1 \!-\! \tau)^2 \|\lambda_k\|^2) B_h^2 \right)$, we obtain:
\begin{align*}
    & \mathbb{E} [ \|x_{k+1} - x^* \|^2] \\
    & \le {(1 - \mu \alpha_k)} \mathbb{E} [\|x_k - x^*\|^2] \!-\! 2 \alpha_k \mathbb{E} [ ( \mathcal{L}_{\rho,\tau} (x_{k+1}; \lambda_{k}) - \mathcal{L}_{\rho,\tau} (x^*; \lambda_{k}))] + \mathbb{E}[\alpha_k^2 B_{\tau, \lambda_k}^2]\\
    & = {(1 - \mu \alpha_k)} \mathbb{E} [\|x_k - x^*\|^2] - 2\alpha_k \mathbb{E}\left[F(x_{k+1}) - F^* + \Psi_{\rho, \tau} (h(x_{k+1}); \lambda_k ) - \Psi_{\rho, \tau} (h(x^*); \lambda_k ) \right] + \mathbb{E}[\alpha_k^2 B_{\tau, \lambda_k}^2] \\
    & = {(1 - \mu \alpha_k)} \mathbb{E} [\|x_k - x^*\|^2] + 2\alpha_k \mathbb{E} \left[\Psi_{\rho, \tau} (h(x^*); \lambda_k ) \right] - 2\alpha_k \mathbb{E}\left[F(x_{k+1}) - F^* +\frac{1 - \tau}{m} \sum_{j=1}^m \lambda^j h_j (x_{k+1})\right]\\
    & \quad + 2\alpha_k \mathbb{E}\left[\frac{1 - \tau}{m} \sum_{j=1}^m \lambda^j h_j (x_{k+1}) - \Psi_{\rho, \tau} (h(x_{k+1}); \lambda_k )\right] + \mathbb{E}[\alpha_k^2 B_{\tau, \lambda_k}^2].
\end{align*}
After using the result from Lemma \ref{lemma_dualupdate}, we get:
\begin{align*}
    & \mathbb{E} [ \|x_{k+1} - x^* \|^2 ] \\
    &   \le (1 - \mu \alpha_k) \mathbb{E}[\|x_k - x^*\|^2]  + 2\alpha_k \mathbb{E} \left[\Psi_{\rho, \tau} (h(x^*); \lambda_k ) \right] + \frac{\alpha_k}{\rho m} \mathbb{E}[\|\lambda_k - \lambda\|^2] + \mathbb{E}[\alpha_k^2 B_{\tau, \lambda_k}^2]    \\
    & \quad  {- \frac{\alpha_k}{\rho m} \mathbb{E}[\|\lambda_{k+1} - \lambda\|^2] - 2\alpha_k \mathbb{E} \left[F(x_{k+1}) - F^* +\frac{1 - \tau}{m} \sum_{j=1}^m \lambda^j h_j (x_{k+1})\right] + 2\alpha_k \frac{\tau}{\rho m}\mathbb{E}[\|\lambda\|^2] .}
\end{align*}
Using the result from Lemma \ref{lem:negative} by noticing that  $ {\lambda}_{k} \ge 0$ (see Lemma \ref{lem:nonnegDual}) and $h_j(x^*) \le 0$, for all $j = 1 :m$, we have:
\begin{align*}
    & {\mathbb{E} [ \|x_{k+1} - x^* \|^2 ] + \frac{\alpha_k}{\rho m} \mathbb{E}[\|\lambda_{k+1} - \lambda\|^2] + 2\alpha_k \mathbb{E} \left[F(x_{k+1}) - F^* +\frac{1 - \tau}{m} \sum_{j=1}^m \lambda^j h_j (x_{k+1})\right]  }   \\
    &  \le (1 - \mu \alpha_k) \mathbb{E}[\|x_k - x^*\|^2] + \frac{\alpha_k}{\rho m} \mathbb{E}[\|\lambda_k - \lambda\|^2] + \mathbb{E}[\alpha_k^2 B_{\tau, \lambda_k}^2] + 2\alpha_k \frac{\tau}{\rho m}\mathbb{E}[\|\lambda\|^2].
\end{align*}
Hence, we obtained the required recurrence.  \Halmos
\endproof 

\noindent Notice that the constants $L_{\rho,\tau}^{\lambda_k}$ and $ B_{\tau, \lambda_k}^2$ in Theorem \ref{lemma_Main} depend on the dual iterate $\lambda_k$. This dependency may lead to $L_{\rho,\tau}^{\lambda_k}, B_{\tau, \lambda_k}^2$ being undefined when the dual iterates are unbounded. To address this issue, in the following two subsections, we consider two cases: one where the dual iterates $\lambda_k$ are assumed  bounded, and another, where we prove that the dual iterates $\lambda_k$ are  bounded.


\subsection{Convergence analysis when the dual iterates are assumed to be bounded}
Now, we are ready to provide the convergence rates for the iterates generated  by  SGDPA algorithm when the dual updates are assumed bounded. We state this as an assumption. 

\begin{assumption}\label{Ass:dualbdd}
The dual iterates generated by  SGDPA algorithm are bounded, i.e., there exists a positive constant $\mathcal{B}_D>0$ such that the following is true:
    $$\|\lambda_k\| \le \mathcal{B}_D \quad \forall k \geq 0.$$
\end{assumption}

\noindent Note that this assumption has  been considered in the literature before, see e.g., \cite{CohTeb:22, BolTeb:14}.   Then, the  smoothness constant $L_{\rho,\tau}^{\lambda_k}$ and the constant $B_{\tau, \lambda_k}^2$ from Theorem \ref{lemma_Main} being dependent on $\lambda_k$, this assumption will ensure well-definedness of these two constants. Thus, we have the following form of the smoothness constant:
\begin{align*}
L_{\rho,\tau}^{\lambda_k} & = L_f + \rho B_h^2 + \left(\rho M_h + \frac{(1-\tau)}{m} \|\lambda_k\|_1 \right)L_h \le L_f + \rho B_h^2 + \left(\rho M_h + \frac{(1-\tau)}{\sqrt{m}} \|\lambda_k\| \right)L_h \\
& \overset{\text{Assumption }\ref{Ass:dualbdd}}{\le} L_f + \rho B_h^2 + \left(\rho M_h + \frac{(1-\tau)}{\sqrt{m}} \mathcal{B}_D \right)L_h = L_{\rho, \tau}^{\mathcal{B}_D},
\end{align*}
where the inequality follows from the equivalence between the norms, i.e., $\|\lambda_k\|_1 \le \sqrt{m} \|\lambda_k\|$. Also the constant $B_{\tau, \lambda_k}^2$ takes the form:
\begin{align*}
    B_{\tau, \lambda_k}^2 & = 2 B^2_f + 8 (\rho^2 M_h^2+ (1 - \tau)^2 \|\lambda_k\|^2) B_h^2  \le 2 B^2_f + 8 (\rho^2 M_h^2 + (1 - \tau)^2 \mathcal{B}_D^2) B_h^2 = B_{\tau, \mathcal{B}_D}^2.
\end{align*}
\noindent \red{Note that both constants $L_{\rho, \tau}^{\mathcal{B}_D}$ and $B_{\tau, \mathcal{B}_D}^2$ depend on $1-\tau$, not on $\tau$. This property is essential in deriving optimal convergence rates for our stochastic algorithm. }


\subsubsection{Convergence rates when the objective is convex}\label{Sec:5.1.1}

\noindent Now, we are ready to provide the convergence rates of the iterates generated by  SGDPA in the convex case {(i.e., $\mu = 0$ in Assumption \ref{Ass:StrongConv})}. First, let us define the following average sequence generated by SGDPA algorithm:
\begin{align}
\label{avg_seq}
\hat{x}_{k+1} = \frac{\sum_{t=0}^{k-1} \alpha_t x_{t+1}}{S_k}, \; \; \text{ where } \; S_k = \sum_{t=0}^{k-1} \alpha_t.
\end{align}

\begin{theorem}\label{th:nonstrongconv}
    Let  Assumptions \ref{Ass:KKT}, \ref{Ass:Smooth} and \ref{Ass:dualbdd} hold and, additionally,  the bounds \eqref{eq:Bddness} be valid. Furthermore, let $\rho > 0, \tau \in [0,1)$ and  choose a non-increasing positive stepsize sequence $\alpha_k \in \left(0, \frac{2}{L_{\rho, \tau}^{\lambda_k}}\right]$, for all $k\ge 0$. Then, for any $(x^*, \lambda^*)$ (a primal-dual solution satisfying Assumption \ref{Ass:KKT}), we have the following estimates true for the average sequence $\hat{x}_k$ generated by  SGDPA algorithm with $\lambda_0 =\mathbf{0}$:
\begin{align*}
    & \mathbb{E} [|F(\hat{x}_{k+1}) - F^*|] \le \frac{1}{S_k}\left(\mathbb{E}[\|x_{0} - x^* \|^2] +  \frac{9\alpha_0}{2 \rho m} \mathbb{E}[\|\lambda^*\|^2] +B_{\tau, \mathcal{B}_D}^2 \sum_{t=0}^{k-1}\alpha_t^2 \right) + \frac{9\tau}{\rho m} \|\lambda^*\|^2,  \\
    & \mathbb{E} \left[ \frac{1}{m} \sum_{j=1}^m [h_j(\hat{x}_{k+1})]_+ \right] \le\! \frac{1}{2(1\!-\!\tau)S_k} \!\left(\! \mathbb{E}[\|x_{0} \!-\! x^* \|^2] \!+\! \frac{\alpha_0}{\rho m} \mathbb{E}[\|\mathbf{1} \!+ \!\lambda^*\|^2] + B_{\tau, \mathcal{B}_D}^2\sum_{t=0}^{k-1}\alpha_t^2 \!\right)\\
    & \qquad \qquad\qquad \qquad\qquad \quad+ \frac{\tau}{\rho (1-\tau) m} \|\mathbf{1} + \lambda^*\|^2.
\end{align*}
\end{theorem}

\proof{Proof:} 
Noticing that  $\alpha_{k+1} \le \alpha_k$ and for $F$ convex  we have $\mu = 0$ in  Theorem \ref{lemma_Main},  we get:
\begin{align*}
    & \mathbb{E} [ \|x_{k+1} - x^* \|^2] + \frac{\alpha_{k+1}}{\rho m} \mathbb{E} [\|\lambda_{k+1} - \lambda\|^2] + 2\alpha_k \mathbb{E} \left[ F(x_{k+1}) - F^* +\frac{1-\tau}{ m} \sum_{j=1}^m \lambda^j h_j(x_{k+1}) \right]\\
    & \le \mathbb{E} [ \|x_{k} - x^* \|^2] + \frac{\alpha_{k}}{\rho m} \mathbb{E} [\|\lambda_{k} - \lambda\|^2] + 2\alpha_k \frac{\tau}{\rho m} \mathbb{E}[\|\lambda\|^2] + \alpha_k^2 B_{\tau, \mathcal{B}_D}^2  \qquad \forall \lambda \geq 0.
\end{align*}
Now, summing it from $0$ to $k-1$, we obtain:
\begin{align*}
    & \mathbb{E} [ \|x_{k} - x^* \|^2] + \frac{\alpha_{k}}{\rho m} \mathbb{E} [\|\lambda_{k} - \lambda\|^2]  + 2\sum_{t=0}^{k-1}\alpha_t \mathbb{E} \left[ F(x_{t+1}) - F^* +\frac{1-\tau}{m} \sum_{j=1}^m \lambda^j h_j(x_{t+1}) \right] \\
    & \le  \mathbb{E} [\|x_{0} - x^* \|^2] + \frac{\alpha_0}{\rho m} \mathbb{E} [\|\lambda_{0} - \lambda\|^2]  + 2\frac{\tau}{\rho m} \mathbb{E} [\|\lambda\|^2] \sum_{t = 0}^{k-1} \alpha_t +  B_{\tau, \mathcal{B}_D}^2 \sum_{t=0}^{k-1}\alpha_t^2 .
\end{align*}
After dividing the whole inequality by $2 S_k$, also using the linearity of expectation, the definition of the average sequence and the convexity of the functions $F, h_j$, for all $j = 1:m$, and noticing that $\mathbb{E} [ \|x_{k} - x^* \|^2] + \frac{\alpha_{k}}{\rho m} \mathbb{E} [\|\lambda_{k} - \lambda\|^2] \ge 0$, we get:
\begin{align}\label{eq:main_con1}
    & \mathbb{E} \left[ F(\hat{x}_{k+1}) - F^* +\frac{1-\tau}{m} \sum_{j=1}^m  \lambda^j h_j(\hat{x}_{k+1}) \right] \\
    & \le \frac{1}{2S_k}\left( \mathbb{E}[\|x_{0} - x^* \|^2] + \frac{\alpha_{0}}{\rho m} \mathbb{E}[\|\lambda_{0} - \lambda\|^2] 
 + B_{\tau, \mathcal{B}_D}^2 \sum_{t=0}^{k-1}\alpha_t^2 \right) +  \frac{\tau}{\rho m} \mathbb{E}[\|\lambda\|^2].\nonumber
\end{align}
Since $- (\lambda^*)^j h_j(\hat{x}_{k+1}) \ge - (\lambda^*)^j [h_j(\hat{x}_{k+1})]_+$, we have from \eqref{eq:KKTconv} that
\begin{align}\label{eq:main_con2}
    F(\hat{x}_{k+1}) - F^* \ge - \frac{1-\tau}{m} \sum_{j=1}^m  (\lambda^*)^j [h_j(\hat{x}_{k+1})]_+.
\end{align}
Thus, with the choice of $\lambda$ as $\lambda^j = \begin{cases}
    1+(\lambda^*)^j & \text{ if } h_j(\hat{x}_{k+1}) > 0\\
    0 &  \text{ otherwise}
\end{cases}$, for any $j = 1:m$, which does not depend on $\bar{j}_k$ (hence, satisfying the assumptions of Lemma \ref{lemma_dualupdate}), from \eqref{eq:main_con1} and \eqref{eq:main_con2}, we get for  $\lambda_0 = \mathbf{0}$:
\begin{align*}
    & \mathbb{E} \left[ \frac{1-\tau}{m} \sum_{j=1}^m [h_j(\hat{x}_{k+1})]_+ \right] \le  \frac{1}{2S_k}\left( \mathbb{E}[\|x_{0} - x^* \|^2] +  \frac{\alpha_0}{\rho m} \mathbb{E}[\|\mathbf{1}+ \lambda^*\|^2] +  B_{\tau, \mathcal{B}_D}^2 \sum_{t=0}^{k-1}\alpha_t^2 \right) + \frac{\tau}{\rho m} \|\mathbf{1} + \lambda^*\|^2.
\end{align*}
Hence, we obtain the feasibility constraint violation bound. Next, if we let  in \eqref{eq:main_con1}
$\lambda^j = \begin{cases}
3(\lambda^*)^j & \text{ if }  h_j(\hat{x}_{k+1}) > 0\\
    0 & \text{ otherwise}
\end{cases},$ for any $j = 1:m$, which again does not depend on $\bar{j}_k$,  and combining with \eqref{eq:main_con2}, we  obtain:
\begin{align}\label{eq:main_con3}
    & \mathbb{E} \left[ \frac{1-\tau}{m} \sum_{j=1}^m  (\lambda^*)^j [h_j(\hat{x}_{k+1})]_+ \right]\\
    & \le  \frac{1}{4S_k}\left( \mathbb{E}[\|x_{0} - x^* \|^2] + \frac{\alpha_0}{\rho m}\mathbb{E}[\|\lambda_{0} - 3 \lambda^*\|^2] + B_{\tau, \mathcal{B}_D}^2 \sum_{t=0}^{k-1}\alpha_t^2 \right) + \frac{9\tau}{2\rho m} \|\lambda^*\|^2 \nonumber.
\end{align}
Hence, by \eqref{eq:main_con3} and \eqref{eq:main_con2}, we have:
\begin{align*}
    & 2 \mathbb{E} [(F(\hat{x}_{k+1}) - F^*)_{-}] \le \frac{1}{2S_k}\left(\mathbb{E}[\|x_{0} - x^* \|^2] + \frac{\alpha_0}{\rho m}\mathbb{E}[\|\lambda_{0} - 3 \lambda^*\|^2] + B_{\tau, \mathcal{B}_D}^2 \sum_{t=0}^{k-1}\alpha_t^2 \right) + \frac{9\tau}{\rho m} \|\lambda^*\|^2,
\end{align*}
where $(a)_- = \max ( - a, 0)$. Additionally, from \eqref{eq:main_con1} with $\lambda = \mathbf{0}$, it follows that:
\begin{align*}
    & \mathbb{E} [F(\hat{x}_{k+1}) - F^*] \le  \frac{1}{2S_k}\left( \mathbb{E}[\|x_{0} - x^* \|^2] + \frac{\alpha_0}{\rho m}\mathbb{E}[\|\lambda_{0} \|^2] + B_{\tau, \mathcal{B}_D}^2 \sum_{t=0}^{k-1}\alpha_t^2 \right).
\end{align*}
Since we know that $|a| = a + 2(a)_{-}$, for any real number $a$, we have:
\begin{align*}
    & \mathbb{E} [|F(\hat{x}_{k+1}) - F^*|] = \mathbb{E} [F(\hat{x}_{k+1}) - F^*] + 2 \mathbb{E} [(F(\hat{x}_{k+1}) - F^*)_{-}]\\
    & \le \frac{1}{S_k} \left( \mathbb{E}[\|x_{0} - x^* \|^2] +  \frac{\alpha_0}{2 \rho m} (\mathbb{E}[\|\lambda_{0} - 3 \lambda^*\|^2] + \mathbb{E}[\|\lambda_0\|^2]) + B_{\tau, \mathcal{B}_D}^2 \sum_{t=0}^{k-1}\alpha_t^2 \right) + \frac{9\tau}{\rho m} \|\lambda^*\|^2.
\end{align*}
Hence we obtain the required result.   \Halmos
\endproof 

\medskip 


\noindent Note that Theorem \ref{th:nonstrongconv} holds for any positive choice of the regularization parameter $\rho$. Now, depending on the value of perturbation parameter $\tau$, i.e., $\tau = 0$ (pure augmented Lagrangian) or $\tau \in (0, 1)$ (perturbed augmented Lagrangian) we provide the convergence rates of SGDPA when the objective $F$ is convex. First, consider the case when $\tau  = 0$.
\begin{corollary}\label{corSGDPA1}
    Under the same conditions as in Theorem \ref{th:nonstrongconv}, with  $\tau = 0$ and the stepsize $\alpha_k=\frac{\alpha_0}{(k+1)^\gamma}$ for $k \ge 1$, where  $\gamma \in [1/2, 1)$ and \red{$\alpha_0 \in \left(0,\frac{2}{L_{\rho, 0}^{\mathcal{B}_D}} \right]$}, we have the following rates of convergence for $\hat{x}_{k}$ of SGDPA (neglecting the logarithmic terms):
\begin{align*}
    & \mathbb{E} [|F(\hat{x}_{k+1}) - F^*|] \le \mathcal{O} \left(  \frac{1}{(k+1)^{1-\gamma}}\right),  \qquad \mathbb{E} \left[ \frac{1}{m} \sum_{j = 1}^m [h_j (\hat{x}_{k+1} )]_+\right] \le \mathcal{O} \left(  \frac{1}{(k+1)^{1 - \gamma}}\right).
\end{align*}
\end{corollary}

\proof{Proof:} 
For $\tau = 0$ and  $\alpha_k=\frac{\alpha_0}{(k+1)^\gamma}$ for $k \ge 1$, with  $\gamma \in [1/2, 1)$ and $\alpha_0 \in \left(0,\frac{2}{L_{\rho, 0}^{\mathcal{B}_D}} \right]$, we have:
\[  S_k = \sum_{t=1}^{k} \alpha_t \geq {\cal O}(\alpha_0 {(k+1)}^{1-\gamma}) \;\; \text{and} \;\; \sum_{t=1}^{k} \alpha_t^2 \leq 
\begin{cases}
    {\cal O}(\alpha_0^2) \!& \! \text{if } \gamma>1/2 \\
    {\cal O}(\alpha^2_0 \cdot \ln(k+1)) \!& \! \text{if } \gamma=1/2. 
\end{cases} 
\] 
Consequently, for  $\gamma \in [1/2, 1)$,  we obtain from Theorem \ref{th:nonstrongconv} the following rates (neglecting the logarithmic terms):
\begin{align*}
    & \mathbb{E} [|F(\hat{x}_{k+1}) - F^*|] \le \frac{1}{\alpha_0 (k+1)^{1-\gamma}} \left( \mathbb{E}[\|x_{0} - x^* \|^2] + \frac{9 \alpha_0}{2 \rho m} \mathbb{E}[\| \lambda^*\|^2] + B_{0, \mathcal{B}_D}^2 \alpha_0^2 \right),\\
    & \mathbb{E} \left[ \frac{1}{m} \sum_{j = 1}^m [h_j (\hat{x}_{k+1} )]_+\right] \le \frac{1}{2\alpha_0 (k+1)^{1-\gamma}} \left( \mathbb{E}[\|x_{0} - x^* \|^2] + \frac{\alpha_0}{\rho m} \mathbb{E}[\|\mathbf{1} + \lambda^*\|^2] + B_{0, \mathcal{B}_D}^2 \alpha_0^2 \right).
\end{align*}
Hence, we obtain the claimed results.   \Halmos
\endproof 

\medskip

\noindent Next, when the perturbation parameter $\tau \in (0, 1)$, in order to obtain convergence rates for the iterates generated by  SGDPA algorithm we  fix the number of iterations, e.g., to  some positive integer  $K$.  Moreover,  we choose $\tau$ of order $\mathcal{O} \left(\frac{1}{(K+1)^{1/2}}\right)$. Consequently, in this case, the definition of the  average sequence is as follows:
$$\hat{x}_{K+1} = \frac{\sum_{t=0}^{K-1} \alpha_t x_{t+1}}{S_K}, 
\;\;\text{where} \; S_K = \sum_{t=0}^{K-1} \alpha_t.$$

\begin{corollary}\label{corSGDPA2}
Under the same conditions as in Theorem \ref{th:nonstrongconv}, with  $\tau \in (0, 1)$ of order $\mathcal{O} \left(\frac{1}{(K+1)^{1-\gamma}}\right)$ and the stepsize $\alpha_k=\frac{\alpha_0}{(k+1)^\gamma}$ for $1 \le k \le K$, with  $\gamma \in [1/2, 1)$,  and $\alpha_0 \in \left(0,\frac{2}{L_{\rho, \tau}^{\mathcal{B}_D}} \right]$, we have the following rates of convergence for $\hat{x}_{K+1}$ of SGDPA (keeping the dominant terms only and neglecting the logarithmic terms):
\begin{align*}
    & \mathbb{E} [|F(\hat{x}_{K+1}) - F^*|] \le \mathcal{O} \left(  \frac{1}{(K+1)^{1-\gamma}}\right), \qquad \mathbb{E} \left[ \frac{1}{m} \sum_{j = 1}^m [h_j (\hat{x}_{K+1} )]_+\right] \le \mathcal{O} \left(  \frac{1}{(K+1)^{1-\gamma}}\right).
\end{align*}
\end{corollary}

\proof{Proof:} 
Note that for $\tau \in (0, 1)$ and the stepsize $\alpha_k=\frac{\alpha_0}{(k+1)^\gamma}$ for $1 \le k \le K$, with  $\gamma \in [1/2, 1)$ and $\alpha_0 \in \left(0,\frac{2}{L_{\rho, \tau}^{\mathcal{B}_D}} \right]$ (recall that $L_{\rho, \tau}^{\mathcal{B}_D}$ and $B_{\tau, \mathcal{B}_D}^2$ depend on $1-\tau$, not on $\tau$), we have:
\[  S_K = \sum_{t=1}^{K} \alpha_t \geq {\cal O}(\alpha_0 {(K+1)}^{1-\gamma}) \;\; \text{and} \;\; \sum_{t=1}^{K} \alpha_t^2 \leq 
\begin{cases}
    {\cal O}(\alpha_0^2) &  \text{if } \gamma>1/2 \\
    {\cal O}(\alpha^2_0 \cdot \ln(K+1)) \!&  \text{if } \gamma=1/2. 
\end{cases} 
\] 
Consequently, for  $\gamma \in (1/2, 1)$, we obtain from Theorem \ref{th:nonstrongconv} the following rates (keeping only the dominant terms):
\begin{align*}
    & \mathbb{E} [|F(\hat{x}_{K+1}) - F^*|] \le \frac{1}{\alpha_0 (K+1)^{1-\gamma}} \left( \mathbb{E}[\|x_{0} - x^* \|^2] + \frac{9\alpha_0}{2\rho m} \mathbb{E}[\|\lambda^*\|^2] + B_{\tau, \mathcal{B}_D}^2 \alpha_0^2 \right) + \frac{9\tau}{\rho m} \|\lambda^*\|^2,\\
    & \mathbb{E} \left[ \frac{1}{m} \sum_{j = 1}^m [h_j (\hat{x}_{K+1} )]_+\right] \le \frac{1}{2(1 - \tau)\alpha_0 (K+1)^{1-\gamma}} \left( \mathbb{E}[\|x_{0} - x^* \|^2] + \frac{\alpha_0}{\rho m} \mathbb{E}[\|\mathbf{1} + \lambda^*\|^2] + B_{\tau, \mathcal{B}_D}^2 \alpha^2_0 \right) \\
    & \qquad\qquad\qquad\qquad\qquad\qquad+  \frac{\tau}{\rho (1 - \tau) m} \|\mathbf{1} + \lambda^*\|^2,
\end{align*}
and keeping only the dominant terms and the fact that $\tau = \mathcal{O} \left(\frac{1}{(K+1)^{1-\gamma}}\right)$, the results follow.  For the particular choice $\gamma=1/2$, if we choose $\tau = \mathcal{O} \left(\frac{1}{(K+1)^{1/2}}\right)$ and neglect the logarithmic terms, we obtain the following  convergence rates:
\begin{align*}
    & \mathbb{E} [|F(\hat{x}_{K+1}) - F^*|] \le \mathcal{O} \left(  \frac{1}{(K+1)^{1/2}}\right),\qquad \mathbb{E} \left[ \frac{1}{m} \sum_{j = 1}^m [h_j (\hat{x}_{K+1} )]_+\right] \le \mathcal{O} \left(  \frac{1}{(K+1)^{1/2}}\right).
\end{align*}
Thus, we get the results.   \Halmos
\endproof 

\medskip 

\noindent Thus,  both Corollaries \ref{corSGDPA1} and \ref{corSGDPA2} show that  SGDPA algorithm can achieve convergence rate of order $\mathcal{O} (1/k^{1/2})$, when the dual iterates, $\lambda_k$, are assumed to be bounded, which matches the optimal rate of convergence for stochastic gradient based methods  from the literature when the objective is assumed to be convex, see e.g.,  \cite{Xu:20, NecSin:22, Ned:11, SinNec:23, SinNec:24COAP}. Moreover, for any given accuracy $\epsilon$,  SGDPA algorithm can find in this case a stochastic $\epsilon$-approximate optimal solution after $\mathcal{O} \left(\frac{1}{\epsilon^{2}}\right)$ iterations,  provided that  the perturbation parameter $\tau$ is either $0$ or proportional to  $\mathcal{O}(\epsilon) \in (0,1)$.


\subsubsection{Convergence rates when the objective is strongly convex}
\label{Sec:5.1.2}

\noindent Here we consider the case when $\mu > 0$ in Assumption \ref{Ass:StrongConv}, i.e., $F$ is a strongly convex function. 

\begin{theorem}\label{th:strongconv}
    Let Assumptions \ref{Ass:KKT} - \ref{Ass:StrongConv} and \ref{Ass:dualbdd} hold and, additionally, the bounds from \eqref{eq:Bddness} be valid.   Furthermore, let $\rho > 0, \tau \in [0,1)$ and  choose a non-increasing positive stepsize sequence $\red{\alpha_k = \min\left(\frac{2}{L_{\rho, \tau}^{\mathcal{B}_D}}, \frac{2}{\mu(k+1)}\right)}$, for all $k \ge 0$. Further, define $k_0 = \lfloor \frac{L_{\rho, \tau}^{\mathcal{B}_D}}{\mu} - 1 \rfloor$. Then, for any $(x^*, \lambda^*)$ (a primal-dual solution satisfying Assumption \ref{Ass:KKT}), we have the following estimates true for the iterates  generated by  SGDPA algorithm with $\lambda_0 =\mathbf{0}$ when $0\le k \le k_0$:
\begin{align*}
    & \mathbb{E} [ \|x_{k_0+1} - x^* \|^2] + \frac{2}{L_{\rho, \tau}^{\mathcal{B}_D} \rho m} \mathbb{E} [\|\lambda_{k_0+1} - \lambda^*\|^2] \\
    & \le \left(1 - \frac{2 \mu}{L_{\rho, \tau}^{\mathcal{B}_D}} \right)^{k_0 + 1} \mathbb{E} [ \|x_{0} \!-\! x^* \|^2] \!+\! \frac{2 \left( 1 + 2\tau (k_0 + 1) \right)}{\rho m L_{\rho, \tau}^{\mathcal{B}_D} }  \mathbb{E} [\|\lambda^* \|^2] + \left(\frac{2 B_{\tau, \mathcal{B}_D }}{L_{\rho, \tau}^{\mathcal{B}_D}} \right)^2   (k_0+1),
    \end{align*}
    and  for the average sequence $\hat{x}_{k+1} =\frac{\sum_{t=k_0+1}^{k} x_{t+1}}{k - k_0}$ when $k > k_0$: 
    \begin{align*}
    & \mathbb{E} [|F(\hat{x}_{k+1}) - F^*|] \le \frac{1}{k-k_0}\left(\frac{\mu (k_0 +1)}{2} \mathbb{E}[\|x_{k_0+1} - x^*\|^2] + \frac{1}{2 \rho m} \left(\mathbb{E}[\|\lambda_{k_0 + 1} - 3 \lambda^*\|^2] + \mathbb{E} [\|\lambda_{k_0+1}\|^2]\right) \right)\\
    & \quad +  \frac{2 B_{\tau, \mathcal{B}_D }^2}{\mu}\cdot \frac{\mathcal{O}( \ln{(k+1)} - \ln{(k_0+1)})}{k-k_0}  + \frac{9 \tau}{\rho m} \mathbb{E}[\|\lambda^*\|^2], \\ 
    & \mathbb{E} \left[ \frac{1}{m} \sum_{j=1}^m [h_j(\hat{x}_{k+1})]_+ \right] \le \frac{1}{(k - k_0) (1 - \tau)}\left( \frac{\mu (k_0 +1)}{4} \mathbb{E} [\|x_{k_0+1} - x^* \|^2] + \frac{1}{2 \rho m} \mathbb{E} [\|\lambda_{k_0+1} - \lambda^* - \mathbf{1}\|^2]  \right)\\
    & \quad +  \frac{B_{\tau, \mathcal{B}_D }^2}{\mu (1 - \tau)}\cdot \frac{\mathcal{O}( \ln{(k+1)} - \ln{(k_0+1)})}{k - k_0} + \frac{\tau}{\rho m (1 - \tau)} \mathbb{E}[\|\lambda^* \!+\! \mathbf{1}\|^2].
\end{align*}
\end{theorem}
    
\proof{Proof:}
For $k\le k_0$, we have $\alpha_k = \frac{2}{L_{\rho, \tau}^{\mathcal{B}_D}}$. Setting $\lambda = \lambda^*$ in  relation \eqref{eq:KKTconv} from Theorem \ref{lemma_Main}, we obtain:
\begin{align*}
    & \mathbb{E} [ \|x_{k+1} - x^* \|^2] + \frac{2}{\rho m L_{\rho, \tau}^{\mathcal{B}_D} } \mathbb{E} [\|\lambda_{k+1} - \lambda^*\|^2]\\
    & \le \left(1 - \frac{2 \mu}{L_{\rho, \tau}^{\mathcal{B}_D}} \right)\mathbb{E} [ \|x_{k} - x^* \|^2] + \frac{2}{\rho m L_{\rho, \tau}^{\mathcal{B}_D} }\mathbb{E} [\|\lambda_{k} - \lambda^*\|^2] + \frac{4}{L_{\rho, \tau}^{\mathcal{B}_D}}\left( \frac{B_{\tau, \mathcal{B}_D }^2}{L_{\rho, \tau}^{\mathcal{B}_D}}  + \frac{\tau}{\rho m} \mathbb{E}[\|\lambda^*\|^2]\right).
\end{align*}
Summing it from $0$ to $k_0$, we get:
\begin{align*}
    & \mathbb{E} [ \|x_{k_0+1} - x^* \|^2] + \frac{2}{\rho m L_{\rho, \tau}^{\mathcal{B}_D} } \mathbb{E} [\|\lambda_{k_0+1} - \lambda^*\|^2] \\
    & \le \left(1 - \frac{2 \mu}{L_{\rho, \tau}^{\mathcal{B}_D}} \right)^{k_0 + 1} \mathbb{E} [ \|x_{0} \!-\! x^* \|^2] \!+\! \frac{2}{\rho m L_{\rho, \tau}^{\mathcal{B}_D} } \mathbb{E} [\|\lambda_{0} \!-\! \lambda^* \|^2]  + \frac{4}{L_{\rho, \tau}^{\mathcal{B}_D}}\left( \frac{B_{\tau, \mathcal{B}_D }^2}{L_{\rho, \tau}^{\mathcal{B}_D}} + \frac{\tau}{\rho m} \mathbb{E}[\|\lambda^*\|^2]\right) (k_0+1).
\end{align*}
Noticing that $\lambda_0 = \mathbf{0}$, after rearranging the terms, we obtain the first claim. Further, for $k> k_0$, we have $\alpha_k = \frac{2}{\mu (k+1)}$ and thus from Theorem \ref{lemma_Main} we get:
\begin{align*}
    & \mathbb{E} [ \|x_{k+1} - x^* \|^2] + \frac{2}{\mu \rho m (k+1)} \mathbb{E} [\|\lambda_{k+1} - \lambda\|^2] + \frac{4}{\mu (k+1)} \mathbb{E} \left[ F(x_{k+1}) - F^* +\frac{1-\tau}{m} \sum_{j=1}^m \lambda^j h_j(x_{k+1}) \right] \\
    & \le  \left( \frac{k-1}{k+1}\right)\mathbb{E} [\|x_{k} - x^* \|^2] + \frac{2}{\mu \rho m (k+1)} \mathbb{E} [\|\lambda_{k} - \lambda\|^2] + \frac{4}{\mu (k+1)} \left( \frac{B_{\tau, \mathcal{B}_D }^2}{\mu (k+1)}  + \frac{\tau}{ \rho m} \mathbb{E}[\|\lambda\|^2]\right).
\end{align*}
Now, multiply the whole inequality by $(k+1)$ and using the fact that $k-1 \le k$, we get:
\begin{align*}
    & (k+1)\mathbb{E} [ \|x_{k+1} - x^* \|^2] + \frac{2}{\mu \rho m} \mathbb{E} [\|\lambda_{k+1} - \lambda\|^2] + \frac{4}{\mu} \mathbb{E} \left[ F(x_{k+1}) - F^* +\frac{1-\tau}{m} \sum_{j=1}^m \lambda^j h_j(x_{k+1}) \right] \\
    & \le k \mathbb{E} [\|x_{k} - x^* \|^2] + \frac{2}{\mu \rho m} \mathbb{E} [\|\lambda_{k} - \lambda\|^2] + \frac{4}{\mu} \left( \frac{B_{\tau, \mathcal{B}_D }^2}{\mu (k+1)}  + \frac{\tau}{ \rho m} \mathbb{E}[\|\lambda\|^2]\right).
\end{align*}
Summing it from $k_0+1$ to $k$, we obtain:
\begin{align*}
    & (k+1)\mathbb{E} [ \|x_{k+1} - x^* \|^2] + \frac{2}{\mu \rho m } \mathbb{E} [\|\lambda_{k+1} - \lambda\|^2] + \frac{4}{\mu}\sum_{t = k_0+1}^k \mathbb{E} \left[ F(x_{t+1}) - F^* +\frac{1-\tau}{m} \sum_{j=1}^m \lambda^j h_j(x_{t+1}) \right] \\
    & \le (k_0 +1) \mathbb{E} [\|x_{k_0+1} - x^* \|^2] + \frac{2}{\mu \rho m} \mathbb{E} [\|\lambda_{k_0+1} - \lambda\|^2]  + \frac{4 B_{\tau, \mathcal{B}_D }^2}{\mu^2}\sum_{t = k_0+1}^k \frac{1}{(t+1)} + \frac{4\tau}{\mu \rho m} \mathbb{E}[\|\lambda\|^2] (k - k_0).
\end{align*}
After multiplying the whole inequality by $\frac{\mu}{4(k - k_0)}$, using the definition of average sequence ${\hat x}_{k+1}$ and the convexity of the functions $F, h_j$, for all $j = 1:m$, and also  the fact that $\sum_{t = k_0+1}^k \frac{1}{(t+1)} \le \mathcal{O} (\ln{(k+1)} - \ln{(k_0+1)})$, we get:
\begin{align*}   
    & \mathbb{E} \left[ F(\hat{x}_{k+1}) - F^* +\frac{1-\tau}{m} \sum_{j=1}^m  \lambda^j h_j(\hat{x}_{k+1}) \right] \\
    & \le \frac{1}{k - k_0}\left( \frac{\mu (k_0 +1)}{4} \mathbb{E} [\|x_{k_0+1} - x^* \|^2] + \frac{1}{2 \rho m} \mathbb{E} [\|\lambda_{k_0+1} - \lambda\|^2] + \frac{B_{\tau, \mathcal{B}_D }^2}{\mu} \cdot \mathcal{O} \left( \ln \frac{k+1}{k_0+1} \right) \right)   + \frac{\tau}{\rho m} \mathbb{E}[\|\lambda\|^2].\nonumber
\end{align*}
Now, following the same reasoning as in the proof of Theorem \ref{th:nonstrongconv} (after equation \eqref{eq:main_con1}), we get the second claim. \Halmos
\endproof 

\medskip 

\noindent Following a similar discussion as  in the convex case, depending on the value of $\tau$, i.e., either $\tau = 0$ or $\tau \in (0, 1)$, we provide below the convergence rates of SGDPA for the strongly convex case. First, consider the case when $\tau =0$.
\begin{corollary}\label{corSGDPA3}
Under the same conditions as in Theorem \ref{th:strongconv} with $\tau =0$,  we have the following rates of convergence for  the average sequence $\hat{x}_{k+1} =\frac{\sum_{t=k_0+1}^{k} x_{t+1}}{k - k_0}$, for any $k > k_0$, generated by  SGDPA (neglecting the logarithmic terms):
\begin{align*}
    & \mathbb{E} [|F(\hat{x}_{k+1}) - F^*|] \le \mathcal{O} \left(  \frac{1}{k - k_0} \right),\qquad \mathbb{E} \left[ \frac{1}{m} \sum_{j = 1}^m [h_j (\hat{x}_{k+1} )]_+\right] \le \mathcal{O} \left(  \frac{1}{k - k_0} \right).
\end{align*}
\end{corollary}

\proof{Proof:}
For $k> k_0$, the proof will follow immediately from Theorem \ref{th:strongconv}.
\Halmos
\endproof

\medskip 

\noindent Next, when $\tau \in (0, 1)$, to obtain convergence rates for SGDPA algorithm we  fix the number of iterations, e.g., to some positive integer  $K > k_0$. Additionally, we choose $\tau$ of order $\mathcal{O} \left(\frac{1}{K}\right)$. 

\begin{corollary} \label{corSGDPA4}
Under the same conditions as in Theorem \ref{th:strongconv}, with  $\tau \in (0, 1)$ of order $\mathcal{O} \left(\frac{1}{K}\right)$, for some fixed $K > k_0$, we have the following rates for $\hat{x}_{K+1} = \frac{\sum_{t=k_0+1}^{K} x_{t+1}}{K - k_0}$ of SGDPA (keeping the dominant terms only and neglecting the logarithmic terms):
\begin{align*}
    & \mathbb{E} [|F(\hat{x}_{K+1}) - F^*|] \le \mathcal{O} \left(  \frac{1}{K - k_0} \right),\qquad \mathbb{E} \left[ \frac{1}{m} \sum_{j = 1}^m [h_j (\hat{x}_{K+1} )]_+\right] \le \mathcal{O} \left( \frac{1}{K - k_0} \right).
\end{align*}
\end{corollary}

\proof{Proof:}
Note that choosing  $\tau = \mathcal{O} (1/K)$, with  $K > k_0$  fixed and recalling  that $L_{\rho, \tau}^{\mathcal{B}_D}$ and $B_{\tau, \mathcal{B}_D}^2$ depend on $1-\tau$, not on $\tau$,  one can easily see from Theorem \ref{th:strongconv} that the claim follows.
\Halmos
\endproof

\medskip 

\noindent Thus,  from Corollaries \ref{corSGDPA3} and \ref{corSGDPA4}, it is clear that SGDPA algorithm achieves the  rate of convergence $\mathcal{O} (1/k)$ (neglecting logarithmic terms),  which matches the optimal rates of convergence from the literature when the objective function is assumed to be strongly convex, see e.g., \cite{Xu:20, NecSin:22, Ned:11, SinNec:23, SinNec:24COAP}. Notably, for any given $\epsilon$, SGDPA algorithm can find a stochastic $\epsilon$-approximately optimal solution after $\mathcal{O} \left(\frac{1}{\epsilon}\right)$ iterations when the dual iterates are assumed to be bounded and the pertubartion parameter $\tau$ is either $0$ or proportional to  $\mathcal{O}(\epsilon) \in (0,1)$.


\subsection{Convergence analysis when the dual iterates  are proved to be bounded}\label{ch04:S2sss}
Proving boundedness of the dual iterates is not a trivial task in optimization, as pointed out also in   \cite{HalTeb:23}. However, using the perturbation introduced in the dual updates \eqref{eq:dual_update} of SGDPA  we prove below that our iterates $\lambda_k$ are indeed bounded, although the computed bound will have negative effect on the convergence rates of SGDPA, being worse than the rates derived in the previous section.   The following lemma derives a bound on the dual iterates generated by  SGDPA algorithm.
\begin{lemma}\label{lem:dualbound}
    Let us assume that the bounds from \eqref{eq:Bddness} hold. Additionally, let $\rho >0, \tau \in (0, 1)$ and  $\lambda_0 = \mathbf{0}$. Then,  the following is true for the dual iterates generated by SGDPA:
    $$|\lambda_k^{j}| \le \frac{\rho M_h}{\tau}   \quad  \forall j = 1:m,\;\; k \ge 0.$$
\end{lemma}
\proof{Proof:}
We will prove this bound using the induction. First, for $k = 0$, it is easy to see that the required result is true according to our initialization choice. Now, let us assume that the statement is true for some $k>0$, i.e.: 
\[ |\lambda_{k}^{j}| \le \frac{\rho M_h}{\tau} \quad \forall j = 1:m. \]
Next, at the $k^{th}$ iteration, let $ \bar{j}_k \in [m]$ be a uniformly at random chosen index, then it follows from the dual update \eqref{eq:dual_update} of SGDPA that:
\begin{align*}
    |\lambda_{k+1}^{j}| = 
    \begin{cases}
        |\lambda_{k}^{j}| & \text{ when } j \neq \bar{j}_k,\\
        \left| (1 - \tau)\lambda_{k}^{j} + \rho \max \left( -\frac{(1 - \tau) \lambda_{k}^{j}}{\rho}, h_{j}(x_{k+1}) \right)\right| & \text{ when } j = \bar{j}_k.
    \end{cases}
\end{align*}
Since we know that $|\lambda_{k}^{j}| \le \frac{\rho M_h}{\tau}$, for all $j = 1:m$, then we have $|\lambda_{k+1}^{j}| = |\lambda_{k}^{j}| \le \frac{\rho M_h}{\tau}$, when $j \neq \bar{j}_k$. Our proof will be complete if we show that the same bound holds when $j = \bar{j}_k$. Indeed: 
\begin{align}
\label{bound1}
|\lambda_{k+1}^{\bar{j}_k}| & \le |(1 - \tau)\lambda_{k}^{\bar{j}_k}| + \left| \rho \max \left( -\frac{(1 - \tau) \lambda_{k}^{\bar{j}_k}}{\rho}, h_{\bar{j}_k}(x_{k+1}) \right)\right| \overset{\text{Lemma } \ref{lemma_bdd_gradDual}}{\le}  (1 - \tau) |\lambda_{k}^{\bar{j}_k}| + \rho M_h \nonumber \\
& \leq (1-\tau) \frac{\rho M_h}{\tau} + \rho M_h  = \frac{\rho M_h}{\tau},
\end{align} 
where the second  inequality follows from the fact that $x_{k+1} \in \mathcal{Y}$ and also by Lemma \ref{lem:nonnegDual} we have $\lambda_{k}^{\bar{j}_k} \ge 0$, while the last inequality follows from our induction assumption. Hence, we get the desired result.
\Halmos
\endproof

\medskip 

\noindent From the bound of  Lemma \ref{lem:dualbound}, it is clear that one cannot consider the case $\tau = 0$ to bound the dual iterates. As reasoned in the previous subsections, here again using the bound from Lemma \ref{lem:dualbound}, we provide an upper bound for the smoothness constant $L_{\rho,\tau}^{\lambda_k}$ of the gradient of the perturbed augmented Lagrangian \eqref{eq:AugLag}:
\begin{align*}
L_{\rho,\tau}^{\lambda_k} & = L_f + \rho B_h^2 + \left(\rho M_h + \frac{(1-\tau)}{m} \sum_{j =1}^m |\lambda^j_k| \right)L_h \overset{\text{Lemma }\ref{lem:dualbound}}{\le} L_f + \rho B_h^2 + \left(\rho M_h + \frac{(1-\tau)}{m} \sum_{j =1}^m \frac{\rho M_h}{\tau} \right)L_h \\
& = L_f + \rho B_h^2 + \frac{\rho M_h}{\tau} L_h := L_{\rho, \tau}.
\end{align*}
Similarly, an upper bound  on the constant $B_{\tau, \lambda_k}^2$ from Theorem \ref{lemma_Main} will be:
\begin{align*}
    B_{\tau, \lambda_k}^2 & = 2 B^2_f + 8 (\rho^2 M_h^2 + (1 - \tau)^2 \|\lambda_k\|^2) B_h^2  \overset{\text{Lemma }\ref{lem:dualbound}}{\le} 2 B^2_f + 8 \left(\rho^2 M_h^2 + (1 - \tau)^2 \sum_{j =1}^m \frac{\rho^2 M_h^2}{\tau^2}\right) B_h^2 \\
    &  = 2 B^2_f + 8 \rho^2 M_h^2 \left(1 + \frac{(1 - \tau)^2 m}{\tau^2}\right) B_h^2 =  B_{\rho, \tau}^2.
\end{align*}

\noindent \red{In contrast to the bounds $L_{\rho, \tau}^{\mathcal{B}_D}$ and $B_{\tau, \mathcal{B}_D}^2$ from Section 5.1 depending on $1-\tau$, one can notice that now   $L_{\rho, \tau}$ is proportional to $1/\tau$, while $B_{\rho,\tau}^2$ is proportional to $1/\tau^2$. In the next sections we investigate the convergence rates of our scheme under these bounds. }


\subsubsection{Convergence rates when the objective is convex}

\noindent In this case we consider the same average sequence as in \eqref{avg_seq}. Then, we have the following convergence estimates:  
\begin{theorem}\label{th:nonstrongconv2}
Let Assumptions \ref{Ass:KKT}, \ref{Ass:Smooth} and \ref{Ass:dualbdd} hold and, additionally,  the bounds from \eqref{eq:Bddness} be valid. Furthermore, let $\rho > 0, \tau \in (0,1)$ and choose a non-increasing positive stepsize sequence $\alpha_k \in \left(0, \frac{2}{L_{\rho, \tau}^{\lambda_k}}\right]$, for all $k\ge 0$. Then,  we have the following estimates true for the average sequence $\hat{x}_k$ from \eqref{avg_seq} generated by  SGDPA algorithm with initialization $\lambda_0 = \mathbf{0}$:
    \begin{align*}
        & \mathbb{E} [|F(\hat{x}_{k+1}) - F^*|] \le \frac{1}{S_k}\left(\mathbb{E}[\|x_{0} - x^* \|^2] +  \frac{9\alpha_0}{2 \rho m} \mathbb{E}[\|\lambda^*\|^2] +B_{\rho, \tau}^2 \sum_{t=0}^{k-1}\alpha_t^2 \right) + \frac{9\tau}{\rho m} \|\lambda^*\|^2,  \\
        & \mathbb{E} \left[ \frac{1}{m} \sum_{j=1}^m [h_j(\hat{x}_{k+1})]_+ \right] \le\! \frac{1}{2(1\!-\!\tau)S_k} \!\left(\! \mathbb{E}[\|x_{0} \!-\! x^* \|^2] \!+\! \frac{\alpha_0}{\rho m} \mathbb{E}[\|\mathbf{1} \!+ \!\lambda^*\|^2] + B_{\rho, \tau}^2\sum_{t=0}^{k-1}\alpha_t^2 \!\right)\\
        & \qquad \qquad\qquad \qquad\qquad \quad+ \frac{\tau}{\rho (1-\tau) m} \|\mathbf{1} + \lambda^*\|^2.
    \end{align*}
\end{theorem}

\proof{Proof:}
 Noticing that $\alpha_{k+1} \le \alpha_k$ and $\mu = 0$ for $F$ convex, from Theorem \ref{lemma_Main} for any $\lambda \ge 0$ we have the following:
\begin{align*}
    & \mathbb{E} [ \|x_{k+1} - x^* \|^2] + \frac{\alpha_{k+1}}{\rho m} \mathbb{E} [\|\lambda_{k+1} - \lambda\|^2] + 2\alpha_k \mathbb{E} \left[ F(x_{k+1}) - F^* +\frac{1-\tau}{ m} \sum_{j=1}^m \lambda^j h_j(x_{k+1}) \right]\\
    & \le \mathbb{E} [ \|x_{k} - x^* \|^2] + \frac{\alpha_{k}}{\rho m} \mathbb{E} [\|\lambda_{k} - \lambda\|^2] + 2\alpha_k \frac{\tau}{\rho m} \mathbb{E}[\|\lambda\|^2] + \alpha_k^2 B_{\rho, \tau}^2.
\end{align*}
Now, proceeding similarly as in the proof of Theorem \ref{th:nonstrongconv}, we obtain the required result. 
\Halmos
\endproof

\medskip 

\noindent Below, we provide the convergence rates for the average  sequence generated by  SGDPA, when $\tau \in (0, 1)$ (perturbed augmented Lagrangian). To obtain convergence rates  we fix the number of iterations to some positive integer  $K$.  We also  choose $\tau$ of order  $\mathcal{O} \left(\frac{1}{(K+1)^{1/4}}\right)$. Consequently, in this case, the definition of the average sequence is as follows (see \eqref{avg_seq}):
$$\hat{x}_{K+1} = \frac{\sum_{t=0}^{K-1} \alpha_t x_{t+1}}{S_K}, 
\;\;\text{where} \; S_K = \sum_{t=0}^{K-1} \alpha_t.$$
\begin{corollary}
\label{cor5}
    Under the same assumptions of Theorem \ref{th:nonstrongconv2} with  $\tau \in (0, 1)$ of order $\mathcal{O} \left(\frac{1}{(K+1)^{(1-\gamma)/2}}\right)$ and the stepsize $\alpha_k=\frac{\alpha_0}{(k+1)^\gamma}$ for $1 \le k \le K$, with  $\gamma \in [1/2, 1)$ and $\alpha_0 \in \left(0,\frac{2}{L_{\rho, \tau}} \right]$, we have the following rates of convergence for the average sequence $\hat{x}_{K+1}$ generated by SGDPA (keeping the dominant terms only and neglecting the logarithmic terms):
\begin{align*}
    & \mathbb{E} [|F(\hat{x}_{K+1}) - F^*|] \le \mathcal{O} \left(  \frac{1}{(K+1)^{(1 -\gamma)/ 2}}\right),\qquad \mathbb{E} \left[ \frac{1}{m} \sum_{j = 1}^m [h_j (\hat{x}_{K+1} )]_+\right] \le \mathcal{O} \left(  \frac{1}{(K+1)^{(1 -\gamma)/ 2}}\right).
\end{align*}
\end{corollary}
\proof{Proof:}
Note that for $\tau \in (0, 1)$ and the stepsize $\alpha_k=\frac{\alpha_0}{(k+1)^\gamma}$ for $1 \le k \le K$, with  $\gamma \in [1/2, 1)$ and $\alpha_0 \in \left(0,\frac{2}{L_{\rho, \tau}} \right]$, we have:
\[  S_K \!=\! \sum_{t=1}^{K} \alpha_t \!\geq\! {\cal O}(\alpha_0 {(K+1)}^{1-\gamma}) \quad \text{and} \quad \sum_{t=1}^{K} \alpha_t^2 \leq 
\begin{cases}
    {\cal O}(\alpha_0^2) \!& \! \text{if } \gamma>1/2 \\
    {\cal O}(\alpha^2_0 \cdot \ln(K+1)) \!& \! \text{if } \gamma=1/2. 
\end{cases} 
\] 
Consequently, for  $\gamma \in (1/2, 1)$, we obtain from Theorem \ref{th:nonstrongconv2} the following estimates:
\begin{align*}
    & \mathbb{E} [|F(\hat{x}_{K+1}) - F^*|] \le \frac{1}{\alpha_0 (K+1)^{1-\gamma}} \left(\mathbb{E}[\|x_{0} - x^* \|^2] + \frac{9\alpha_0}{2\rho m} \mathbb{E}[\|\lambda^*\|^2] + B_{\rho, \tau}^2 \alpha_0^2 \right) + \frac{9\tau}{\rho m} \|\lambda^*\|^2,\\
    & \mathbb{E} \left[ \frac{1}{m} \sum_{j = 1}^m [h_j (\hat{x}_{K+1} )]_+\right] \le \frac{1}{2(1 - \tau)\alpha_0 (K+1)^{1-\gamma}} \left( \mathbb{E}[\|x_{0} - x^* \|^2] + \frac{\alpha_0}{\rho m} \mathbb{E}[\|\mathbf{1} + \lambda^*\|^2] + B_{\rho, \tau}^2 \alpha_0^2 \right)\\
    & \qquad \qquad\qquad \qquad\qquad \quad + \frac{\tau}{\rho (1 - \tau) m} \|\mathbf{1} + \lambda^*\|^2 .
\end{align*}
\noindent  Recall that $L_{\rho, \tau}$ is proportional to $1/\tau$, hence choosing $\alpha_0$  of the form $\alpha_0 = \frac{2}{L_{\rho, \tau}}  \sim \mathcal{O} (\tau)$ and since $ B_{\rho, \tau}^2 $ is proportional to $1/ \tau^2$, the right hand side in the previous estimates are of order:
\[  \frac{1}{(K+1)^{1-\gamma}} +  \frac{1}{\tau (K+1)^{1-\gamma}} + \tau. \]
Therefore, choosing $\tau = \mathcal{O}\left(\frac{1}{(K+1)^{(1-\gamma)/2}}\right)$ and keeping only the dominant terms, we get the required results. Similarly, for the particular choice $\gamma=1/2$, keeping the dominant terms only and neglecting the logarithmic terms, we obtain the following rates of convergence:
\begin{align*}
    & \mathbb{E} [|F(\hat{x}_{K+1}) - F^*|] \le \mathcal{O} \left(  \frac{1}{(K+1)^{1/4}}\right),\qquad  \mathbb{E} \left[ \frac{1}{m} \sum_{j = 1}^m [h_j (\hat{x}_{K+1} )]_+\right] \le \mathcal{O} \left(  \frac{1}{(K+1)^{1/4}}\right).
\end{align*}
Hence we obtain the statement of the corollary.
\Halmos
\endproof

\medskip 

\noindent Thus, in the convex case, for any given $\epsilon$,  SGDPA algorithm can find a stochastic $\epsilon$-approximate optimal solution after $\mathcal{O} \left(\frac{1}{\epsilon^{4}}\right)$ iterations when the dual iterates, $\lambda_k$, are  bounded as in Lemma \ref{lem:dualbound} and the perturbation parameter $\tau $ is chosen of order $\mathcal{O}(\epsilon)$.


\subsubsection{Convergence rates when the objective is strongly convex}
\noindent Here we consider the case when $\mu > 0$, i.e., $F$ is a strongly convex function. 

\begin{theorem}\label{th:strongconv2}
Let Assumptions \ref{Ass:KKT} - \ref{Ass:StrongConv} and \ref{Ass:dualbdd} hold and, additionally, the bounds from \eqref{eq:Bddness} are valid.  Moreover, let $\rho > 0, \tau \in (0,1)$ and choose a non-increasing positive stepsize sequence $\alpha_k = \min\left(\frac{2}{L_{\rho, \tau}}, \frac{2}{\mu(k+1)}\right)$, for all $k \ge 0$. Further, define $k_0 = \lfloor \frac{L_{\rho, \tau}}{\mu} - 1 \rfloor$. Then, for any $(x^*, \lambda^*)$ (a primal-dual solution satisfying Assumption \ref{Ass:KKT}), we have the following estimates true for the iterates  generated by  SGDPA algorithm with $\lambda_0 =\mathbf{0}$ when $0\le k \le k_0$:
\begin{align*}
    & \mathbb{E} [ \|x_{k_0+1} - x^* \|^2] + \frac{2}{L_{\rho, \tau} \rho m} \mathbb{E} [\|\lambda_{k_0+1} - \lambda^*\|^2] \\
    & \le \left(1 - \frac{2 \mu}{L_{\rho, \tau}} \right)^{k_0 + 1} \mathbb{E} [ \|x_{0} - x^* \|^2] + \frac{2 \left( 1 + 2\tau (k_0 + 1) \right)}{\rho m L_{\rho, \tau} }  \mathbb{E} [\|\lambda^* \|^2] + \left(\frac{2 B_{\rho, \tau}}{L_{\rho, \tau}} \right)^2   (k_0+1).
\end{align*}
 and  for the average sequence $\hat{x}_{k+1} =\frac{\sum_{t=k_0+1}^{k} x_{t+1}}{k - k_0}$ when $k > k_0$:
\begin{align*} 
    & \mathbb{E} [|F(\hat{x}_{k+1}) - F^*|] \le \frac{1}{k-k_0}\left(\frac{\mu (k_0 +1)}{2} \mathbb{E}[\|x_{k_0+1} - x^*\|^2] + \frac{1}{2 \rho m} \left(\mathbb{E}[\|\lambda_{k_0 + 1} - 3 \lambda^*\|^2] + \mathbb{E} [\|\lambda_{k_0+1}\|^2]\right) \right)\\
    & \qquad +  \frac{2 B_{\rho, \tau}^2}{\mu}\cdot \frac{\mathcal{O}( \ln{(k+1)} - \ln{(k_0+1)})}{k - k_0}  + \frac{9 \tau}{\rho m} \mathbb{E}[\|\lambda^*\|^2], \\ 
    & \mathbb{E} \left[ \frac{1}{m} \sum_{j=1}^m [h_j(\hat{x}_{k+1})]_+ \right] \le \frac{1}{(k - k_0) (1 - \tau)}\left( \frac{\mu (k_0 +1)}{4} \mathbb{E} [\|x_{k_0+1} - x^* \|^2] + \frac{1}{2 \rho m} \mathbb{E} [\|\lambda_{k_0+1} - \lambda^* - \mathbf{1}\|^2]  \right)\\
    & \qquad +  \frac{B_{\rho, \tau}^2}{\mu (1 - \tau)}\cdot \frac{\mathcal{O}( \ln{(k+1)} - \ln{(k_0+1)})}{k - k_0} + \frac{\tau}{\rho m (1 - \tau)} \mathbb{E}[\|\lambda^* + \mathbf{1}\|^2].
\end{align*}
\end{theorem}
    
\proof{Proof:}
Following the same reasoning as in the proof of  Theorem \ref{th:strongconv} for both cases, i.e., when $k\le k_0$ and when $k>k_0$, we obtain the results.
\Halmos
\endproof

\medskip 

\noindent \red{Recall that $L_{\rho, \tau}$ is proportional to $1/\tau$, hence  $\mu k_0$  is of  order  $1/\tau$;  since $ B_{\rho, \tau}^2 $ is proportional to $1/ \tau^2$, then  $B_{\rho, \tau}^2/\mu$ is of order $k_0/\tau$ and  the right hand side in the  estimates of Theorem  \ref{th:strongconv2} are of order:
\[  \frac{1}{k - k_0} +  \frac{1}{\tau (k-k_0)} +  \frac{k_0(\ln k - \ln k_0)}{\tau (k-k_0)} + \tau. \]}

\noindent Therefore, fixing the number of iteration  $K>k_0$ and considering $\tau$ of order  $\mathcal{O} (1/\sqrt{K})$, we get the following rates of convergence for SGDPA (keeping the dominant terms only and neglecting the logarithmic terms):
\begin{align*}
    & \mathbb{E} [|F(\hat{x}_{K+1}) - F^*|] \le \mathcal{O} \left(  \frac{1}{\sqrt{K}} \right),\quad  \mathbb{E} \left[ \frac{1}{m} \sum_{j = 1}^m [h_j (\hat{x}_{K+1} )]_+\right] \le \mathcal{O} \left(  \frac{1}{\sqrt{K}} \right).
\end{align*}

\noindent In conclusion, for strongly convex objective,  for any given $\epsilon$,  SGDPA algorithm can find a stochastic $\epsilon$-approximate optimal solution after $\mathcal{O} \left(\frac{1}{\epsilon^{2}}\right)$ iterations when the dual iterates, $\lambda_k$, are  bounded as in Lemma \ref{lem:dualbound} and the perturbation parameter $\tau$ is chosen of order   $\mathcal{O}(\epsilon)$. Notice that although our perturbed augmented Lagrangian framework allows us to prove that the dual iterates generated by SGDPA algorithm are bounded, the derived bound depends on the perturbation parameter $\tau$ (see Lemma \ref{lem:dualbound}), which leads to  worse rates of convergence than when the multipliers are assumed to be bounded. It is an open question whether one can derive a better  bound for the dual iterates  than the one obtained in Lemma \ref{lem:dualbound} (e.g., independent on $\tau$).


\section{Numerical experiments }
In this section, we test the performance of our algorithm on convex quadratically constrained quadratic programs (QCQPs) using synthetic and real data. A QCQP problem is given as:
\begin{align}
\label{qcqp}
& \min\limits_{x \in \mathcal{Y} \subseteq \mathbb{R}^n} \; \frac{1}{2}x^TQ_f x + q_f^T x \quad  \text{s.t.} \;\;\;  
	   \frac{1}{2} x^T Q_i x + q_i^T x - b_i \le 0 \;\; \forall i = 1:m. 
\end{align}
We compare our Stochastic Gradient Descent and Perturbed Ascent (SGDPA) algorithm to the  deterministic LALM algorithm \cite{Xu:21}, the stochastic PDSG method (although developed for nonsmooth functional constraints, we adapted it for smooth constraints) \cite{Xu:20} and the  commercial optimization software package, FICO \cite{FICO} (having a specialized solver for QCQPs). The implementation details are conducted using MATLAB R2023b on a laptop equipped with an i5 CPU operating at 2.1 GHz and 16 GB of RAM.


 \subsection{\red{Selection of  initial stepsize}}
 \label{sec_adaptive}
\red{The convergence results from the previous sections  rely on the assumption that the initial stepsize $\alpha_0$ is below a certain threshold, see e.g., Corollaries \ref{corSGDPA1}-\ref{cor5} and Theorem  \ref{th:strongconv2}. However, in practice, determining this threshold beforehand poses challenges as it depends on unknown parameters of the  problem's data and the algorithm's parameters. To overcome this challenge, we propose below an outer algorithm that repeatedly calls Algorithm \ref{algorithmSGD} for a fixed number of iterations, denoted as $K_t$, using an initial stepsize   $\alpha_0^t$. If Algorithm \ref{algorithmSGD} does not yield an $\epsilon$-approximate solution for the problem \eqref{eq:prob} within $K_t$ iterations, then  $K_t$ is increased  and the initial stepsize $\alpha_0^t$ is decreased geometrically. Specifically, we set $K_{t+1} = \zeta_1 K_t$ and $\alpha_0^{t+1} = \zeta_2 \alpha_0^t$, where $\zeta_1 > 1$ and $\zeta_2 < 1$. The resulting procedure can be summarized in the following algorithm:}

\begin{algorithm}[H]
\caption{SGDPA with trial values of $\alpha_0$}
\label{algorithmSGD_1}
\begin{algorithmic}[1]
\State $\text{Choose} \; x_{-1}^* \in \mathcal{Y}, \lambda_{-1}^* \geq 0, \tau \in [0, 1), \; \rho > 0, \; K_0>0, \zeta_1>1, \; \zeta_2 < 1  \; \text{ and initial stepsize } \; \alpha_0^0>0$.
\State   $t \gets 0$
\State \textbf{while} $\epsilon$-approximate solution not found \textbf{do}
\vspace*{-0.3cm}
\begin{align*}
    & \hspace{0.2cm} 3.1\!:\; \text{Call Algorithm \ref{algorithmSGD} with initial stepsize } \;  \alpha_0= \alpha_0^t \; \text{ and warm start } \; (x_0,\lambda_0)=(x_{t-1}^*,\lambda_{t-1}^*)  \\
    & \hspace{1cm}  \text{for } \;  K_t \; \text{ iterations yielding } \; (x_{K_t}, \lambda_{K_t})  \\
    & \hspace{0.2cm} 3.2\!:\;  \text{Update } \; (x_{t}^{*}, \lambda_{t}^{*}) \gets (x_{K_t}, \lambda_{K_t}) \\
    & \hspace{0.2cm} 3.3\!: \; \text{Update } \; K_{t+1} \gets \zeta_1 K_t \; \text{ and } \;  \alpha_0^{t+1} \gets \zeta_2 \alpha_0^t\\
    & \hspace{0.2cm} 3.4\!: \; \text{Update outer iteration } \; t \gets t+1
\end{align*}
\vspace*{-0.8cm}
\State \textbf{end while}
\end{algorithmic}
\end{algorithm}

\vspace*{-0.5cm}

\noindent \red{In our numerical experiments  we considered this adaptive variant of SGDPA (Algorithm \ref{algorithmSGD_1}) and we noticed that in practice it requires, in average, 1-2 restarts (outer iterations). At   each outer iteration $t$,  we choose a stepsize sequence $\alpha_k^t = \frac{\alpha_0^t}{\sqrt{k+1}}$ (convex case) or $\alpha_k^t = \min\left( \alpha_0^t, \frac{2}{\mu(k+1)} \right)$ (strongly convex case) for the inner iterations $k$.}


\subsection{Solving  synthetic convex QCQPs} 
\label{ch04:S2SS5}

\noindent To solve synthetic QCQPs, \eqref{qcqp}, we consider two settings for the objective function, i.e., convex and strongly convex, and convex functional constraints. Also, we consider $\mathcal{Y} = \mathbb{R}^n_+$. To generate the synthetic data we follow the same technique as in \cite{SinNec:24COAP}. 
This means to generate $m$ convex constraints in \eqref{qcqp}, choosing $Q_i = Y_i^T D_i Y_i$, where $Y_i$'s are randomly generated orthogonal $n\times n$ matrices using \texttt{orth} in Matlab and each $D_i$ is a diagonal matrix with $n/10$ zero entries on the diagonal and the rest of the entries are uniformly random generated in $(0, 1)$. Moreover, for the convex objective, the positive semidefinite matrix $Q_f = Y_f^T D_f Y_f$ is generated in the same fashion as  $ Q_i$'s, and for the strong convex case each entry in the diagonal of the matrix $D_f$ is  nonzero uniformly random generated in $(0,1)$. Furthermore, the vectors $q_f, q_i \in \mathbb{R}^n$, for all $i = 1:m$, are generated from uniform distributions. For the scalars $b_i\in \mathbb{R}$ for all $i = 1:m$, we have considered two scenarios: for the first half of the table we generated uniformly at random an initial point $x_0$ and then chose $b_i$'s as $b_i = \frac{1}{2} x_0^T Q_i x_0 + q_i^T x + 0.1$; for the second half of the table we generated $b_i$'s uniformly at random. Note that in both scenarios the QCQPs are feasible, since in the first case $x_0$ is a feasible point, while in the second case, $0$ is always a feasible point. 

\medskip 

\noindent  For SGDPA we  set $\tau$ to either $0$ or to the  accuracy's value and  $\rho = 10$. We stop SGDPA, LALM and PDSG when  $\|\max(0, h(x))\|^2 \le 10^{-2}$ and $|F(x) - F^*| \leq 10^{-2}$ (with $F^*$ computed via FICO solver) or when $\max (\|x_{k+1} - x_k\|^2, \ldots, \|x_{k-M+1} - x_{k-M}\|^2) \le 10^{-3}$, with  $M=10$ (when FICO does not solve the QCQP due to license limitations).

\medskip

\noindent In Table \ref{table:SGDPA}, we provide the average CPU times  in seconds (out of $10$ runs) together with the standard deviation (std)\footnote{\red{Standard deviation is computed as: $ \text{std} = \left(1/10 \sum_{i=1}^{10}  (\text{cpu}(i) - E(\text{cpu}))^2 \right)^{1/2}$.}} for the proposed algorithm SGDPA and for  the  algorithm PDSG, and CPU times  for the  method LALM and for the commercial solver FICO for various dimensions $(n, m)$ in the QCQPs \eqref{qcqp}. The best average time achieved by an algorithm to solve a given QPQC is written in bold. One can observe that SGDPA always outperforms LALM and PDSG, however, FICO sometimes has better behavior than SGDPA (usually on small dimensions). The $"*"$ in LALM column means that the algorithm is not finished in $7$ hours, and in FICO column means that FICO was not solving the problem due to academic license limitation on the problem dimension.

\medskip

\begin{figure}[H]
    \centering
    \includegraphics[height=6cm, width=6.5cm]{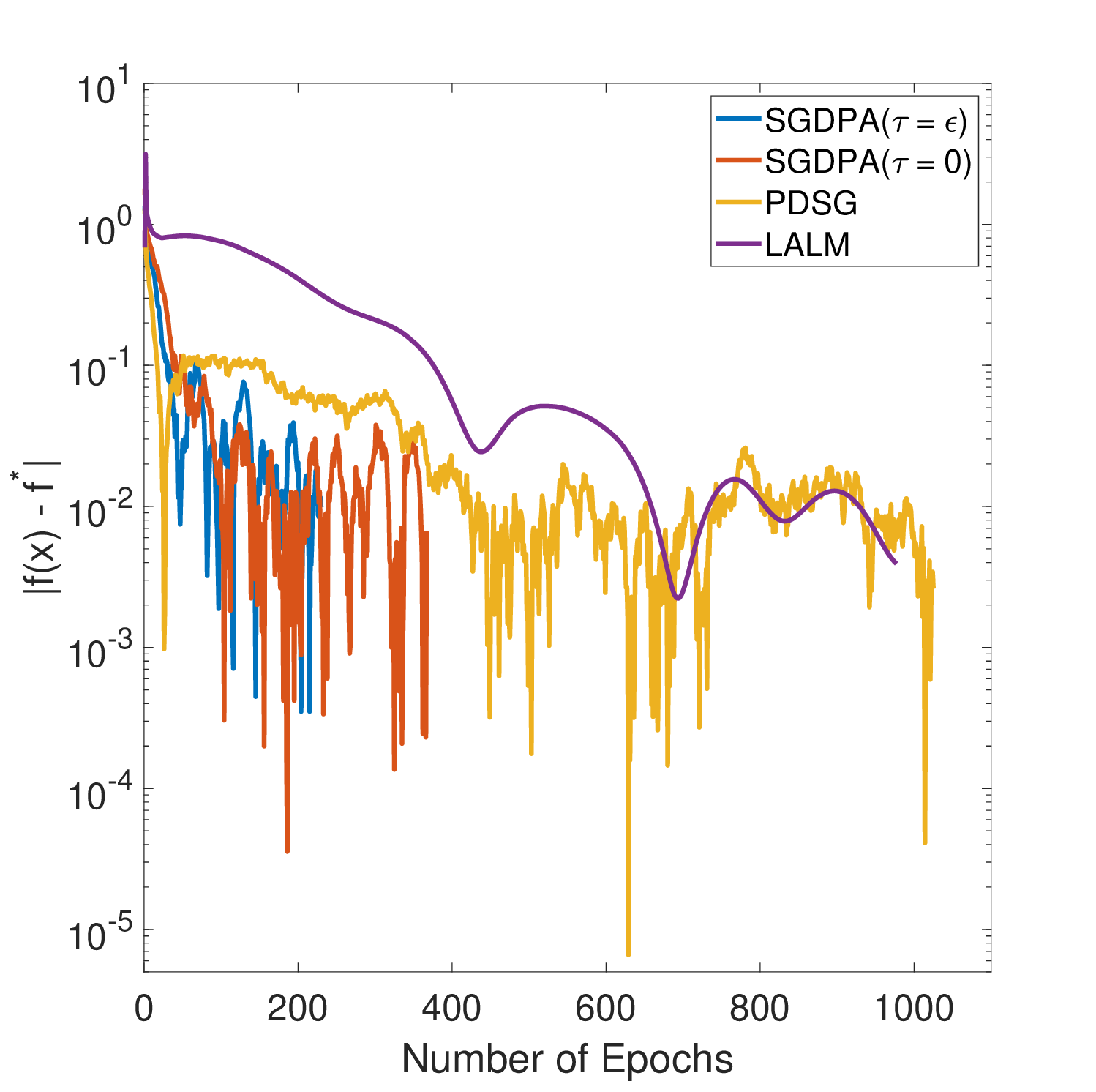}
    \includegraphics[height=6cm, width=6.5cm]{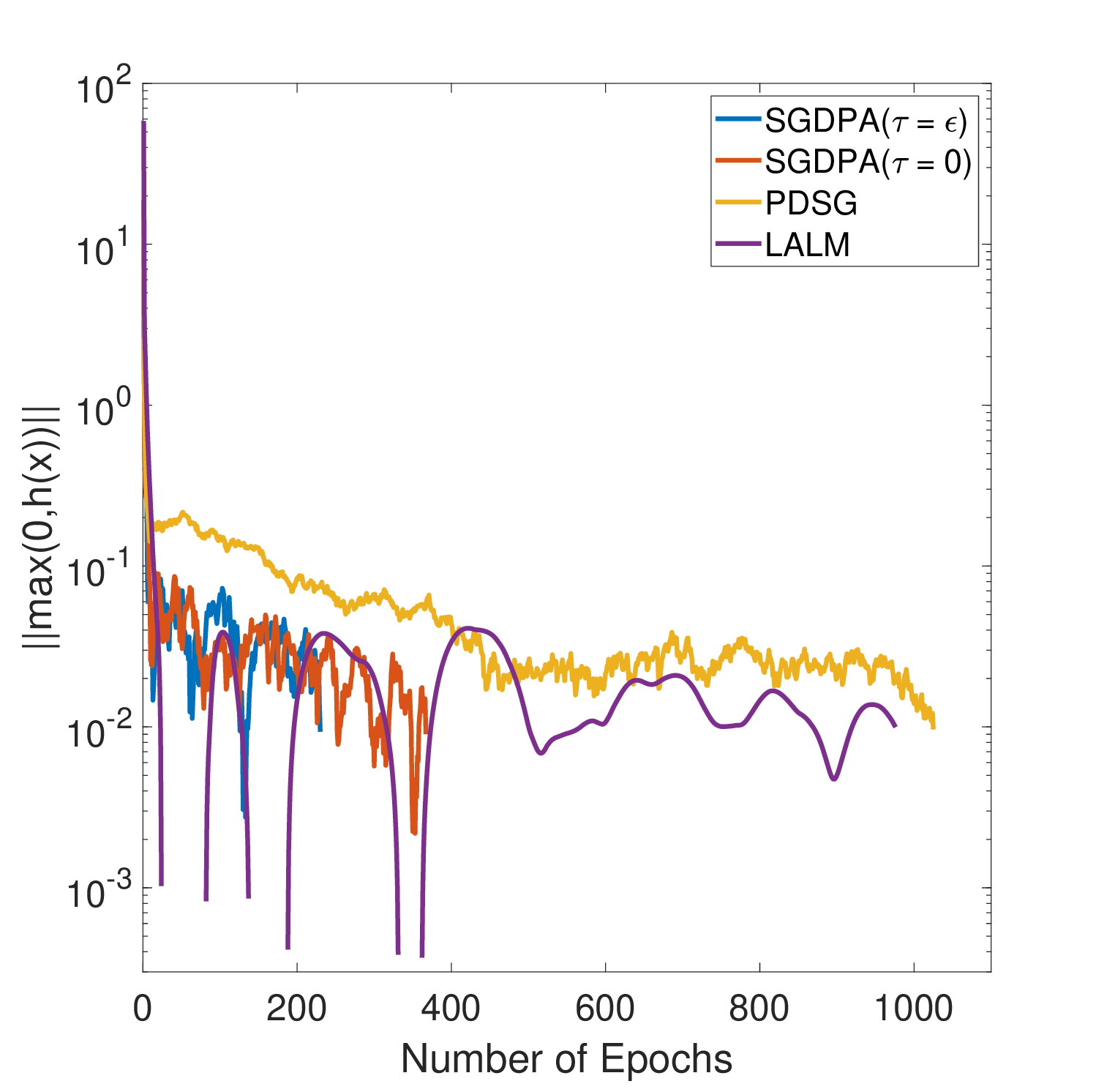}
    \caption{Behaviour of SGDPA  for two choices of $\tau$ ($0$ and $\epsilon=10^{-2}$) and comparison with LALM and PDSG in terms of optimality (left) and feasibility (right) for $n=100,  m = 100$ along the number of epochs.}
    \label{fig:1SGDPA}
\end{figure}

\noindent In Figure \ref{fig:1SGDPA}, we plot the behaviors of SGDPA algorithm (with $\tau = \epsilon = 10^{-2}$ and $\tau = 0$), LALM and PDSG methods along the number of epochs, focusing on optimality (left) and feasibility trends (right) for dimension $n=100$ and $m=100$ of the QCQP having a convex objective. For SGDPA and PDSG, one epoch means $m$ number of iterations (recall that $m$ is the number of constraints); while LALM being a deterministic algorithm,  we consider each iteration as an epoch. It is evident that SGDPA algorithm significantly outperforms LALM and PDSG, demonstrating a clear advantage in terms of optimality and feasibility over the number of epochs. Furthermore, the breaks observed in the feasibility plots occur because all the algorithms are feasible at certain iterations, meaning $\|\max(0, h(x))\| = 0$. However, since the algorithm has not yet satisfied the optimality criterion at those points, it continues. During this process, SMBA, PDSG and LALM can become infeasible again, as these algorithms do not guarantee feasibility at each step.

\begin{table}[H]
\centering
\begin{tabular}{|ccccccccc|}
\hline
\multicolumn{1}{|c|}{\multirow{2}{*}{\textbf{\begin{tabular}[c]{@{}c@{}}Data\\ (n, m)\end{tabular}}}} & \multicolumn{4}{c|}{SGDPA}                                                                                                                                                              & \multicolumn{1}{c|}{\multirow{2}{*}{LALM}} & \multicolumn{2}{c|}{PDSG} & \multicolumn{1}{l|}{\multirow{2}{*}{FICO}} \\ 
\multicolumn{1}{|c|}{}                                                                                & \multicolumn{1}{c|}{$\tau = 0$}       & \multicolumn{1}{c|}{std}                  & \multicolumn{1}{c|}{$\tau = 10^{-2}$}                & \multicolumn{1}{c|}{std}                  & \multicolumn{1}{c|}{}      & \multicolumn{1}{c|}{time}   & \multicolumn{1}{c|}{std}             & \multicolumn{1}{l|}{}                      \\ \hline


\multicolumn{1}{|c|}{$10^2, 10^2$}                       & \multicolumn{1}{c|}{0.68} & \multicolumn{1}{c|}{0.17} & \multicolumn{1}{c|}{\textbf{0.59}} & \multicolumn{1}{c|}{0.11} & \multicolumn{1}{c|}{0.8}   & \multicolumn{1}{c|}{2.59}      & \multicolumn{1}{c|}{0.55}        & \multicolumn{1}{c|}{0.59}               \\ 
\hline



\multicolumn{1}{|c|}{$10^2, 10^3$}                       & \multicolumn{1}{c|}{3.59}                 & \multicolumn{1}{c|}{0.44} & \multicolumn{1}{c|}{\textbf{2.2}}                 & \multicolumn{1}{c|}{0.51} & \multicolumn{1}{c|}{4.69}     & \multicolumn{1}{c|}{196.48}    & \multicolumn{1}{c|}{10.26}       & \multicolumn{1}{c|}{2.29}   \\ 
\hline

\multicolumn{1}{|c|}{$10^2, 5\times 10^3$}                              & \multicolumn{1}{c|}{10.92}        & \multicolumn{1}{c|}{0.53} & \multicolumn{1}{c|}{8.37}                          & \multicolumn{1}{c|}{0.61} & \multicolumn{1}{c|}{16.52}     & \multicolumn{1}{c|}{228.69}   & \multicolumn{1}{c|}{9.24}       & \multicolumn{1}{c|}{\textbf{8.10}}   \\ 
\hline

\multicolumn{1}{|c|}{$10^3, 10^2$}                        & \multicolumn{1}{c|}{13.59}        & \multicolumn{1}{c|}{0.64} & \multicolumn{1}{c|}{\textbf{11.26}}                          & \multicolumn{1}{c|}{0.29}                  & \multicolumn{1}{c|}{21.89}     & \multicolumn{1}{c|}{28.54}    & \multicolumn{1}{c|}{0.89}         & \multicolumn{1}{c|}{20.74}  \\ 
\hline

\multicolumn{1}{|c|}{$10^3, 10^3$}                      & \multicolumn{1}{c|}{110.68}        & \multicolumn{1}{c|}{0.88}                 & \multicolumn{1}{c|}{\textbf{89.62}}                          & \multicolumn{1}{c|}{0.76}                 & \multicolumn{1}{c|}{256.9}       & \multicolumn{1}{c|}{854.51}   & \multicolumn{1}{c|}{5.24}     & \multicolumn{1}{c|}{205.52}\\ 
\hline

\multicolumn{1}{|c|}{$10^3, 5\times10^3$}                & \multicolumn{1}{c|}{\textbf{170.29}}                 & \multicolumn{1}{c|}{0.48}                 & \multicolumn{1}{c|}{188.36}                 & \multicolumn{1}{c|}{0.85}                 & \multicolumn{1}{c|}{$1.05\times 10^3$}    & \multicolumn{1}{c|}{$1.8\times 10^3$}    & \multicolumn{1}{c|}{5.68}       & \multicolumn{1}{c|}{557.44}  \\ 
\hline


\multicolumn{1}{|c|}{$10^3, 10^4$}                       & \multicolumn{1}{c|}{370.12}        & \multicolumn{1}{c|}{1.2}                 & \multicolumn{1}{c|}{\textbf{350.45}}                          & \multicolumn{1}{c|}{1.35}                 & \multicolumn{1}{c|}{*}  & \multicolumn{1}{c|}{$2.6\times 10^3$}   & \multicolumn{1}{c|}{8.19}   & \multicolumn{1}{c|}{*}  \\ 
\hline


\multicolumn{1}{|c|}{$5 \times 10^3, 10^4$}              & \multicolumn{1}{c|}{578.21}        & \multicolumn{1}{c|}{2.1}                 & \multicolumn{1}{c|}{\textbf{521.8}}                          & \multicolumn{1}{c|}{2.5}                 & \multicolumn{1}{c|}{*}  & \multicolumn{1}{c|}{$8.6 \times 10^3$}   & \multicolumn{1}{c|}{8.91}   & \multicolumn{1}{c|}{*}  \\ \hline


\multicolumn{1}{|c|}{$10^4, 10^4$}                      & \multicolumn{1}{c|}{$3.6 \times 10^3$}        & \multicolumn{1}{c|}{7.5}                 & \multicolumn{1}{c|}{\textbf{2.8} $\times$ \textbf{$10^3$}}                          & \multicolumn{1}{c|}{8.3}                 & \multicolumn{1}{c|}{*}  & \multicolumn{1}{c|}{$1.9 \times 10^4$}   & \multicolumn{1}{c|}{9.4}   & \multicolumn{1}{c|}{*}  \\ \hline
\hline


\multicolumn{1}{|c|}{$10^2, 10^2$}                                                                    & \multicolumn{1}{c|}{0.71}        & \multicolumn{1}{c|}{0.29} & \multicolumn{1}{c|}{\textbf{0.63}}            & \multicolumn{1}{c|}{0.2} & \multicolumn{1}{c|}{1.47}      & \multicolumn{1}{c|}{{3.1}}  & \multicolumn{1}{c|}{0.71}          & \multicolumn{1}{c|}{0.65}      \\ 
\hline



\multicolumn{1}{|c|}{$10^2, 10^3$}                     & \multicolumn{1}{c|}{3.26}        & \multicolumn{1}{c|}{0.4}   & \multicolumn{1}{c|}{{2.65}}           & \multicolumn{1}{c|}{0.57} & \multicolumn{1}{c|}{6.32}     & \multicolumn{1}{c|}{254.03}   & \multicolumn{1}{c|}{11.47}          & \multicolumn{1}{c|}{\textbf{2.27}}  \\ \hline

\multicolumn{1}{|c|}{$10^2, 5 \times 10^3$}                              & \multicolumn{1}{c|}{11.48}        & \multicolumn{1}{c|}{0.62} & \multicolumn{1}{c|}{{9.89}}                          & \multicolumn{1}{c|}{0.46}                 & \multicolumn{1}{c|}{21.17}     & \multicolumn{1}{c|}{311.26}    & \multicolumn{1}{c|}{8.93}       & \multicolumn{1}{c|}{\textbf{8.33}}      \\ \hline

\multicolumn{1}{|c|}{$10^3, 10^2$}                        & \multicolumn{1}{c|}{14.23}        & \multicolumn{1}{c|}{0.51}                 & \multicolumn{1}{c|}{\textbf{12.03}}                          & \multicolumn{1}{c|}{0.45}                 & \multicolumn{1}{c|}{35.34}         & \multicolumn{1}{c|}{48.73}   & \multicolumn{1}{c|}{1.16}      & \multicolumn{1}{c|}{21.16}     \\ \hline

\multicolumn{1}{|c|}{$10^3, 10^3$}                      & \multicolumn{1}{c|}{133.12}                 & \multicolumn{1}{c|}{0.97}                 & \multicolumn{1}{c|}{\textbf{108.84}}                 & \multicolumn{1}{c|}{0.89}                 & \multicolumn{1}{c|}{333.62}      & \multicolumn{1}{c|}{$1.07 \times 10^3$}   & \multicolumn{1}{c|}{6.23}       & \multicolumn{1}{c|}{208.18}  \\ \hline


\multicolumn{1}{|c|}{$10^3, 5 \times 10^3$}             & \multicolumn{1}{c|}{238.32}                & \multicolumn{1}{c|}{1.02}                 & \multicolumn{1}{c|}{\textbf{202.45}}                & \multicolumn{1}{c|}{0.94}                 & \multicolumn{1}{c|}{$1.1 \times 10^3$}   & \multicolumn{1}{c|}{$2.1 \times 10^3$}   & \multicolumn{1}{c|}{7.1}    & \multicolumn{1}{c|}{861.70}     \\  \hline

\multicolumn{1}{|c|}{$10^3, 10^4$}                      & \multicolumn{1}{c|}{410.45}               & \multicolumn{1}{c|}{1.32}                 & \multicolumn{1}{c|}{\textbf{390.23}}               & \multicolumn{1}{c|}{1.28}                 & \multicolumn{1}{c|}{*}  & \multicolumn{1}{c|}{$4.3 \times 10^3$}   & \multicolumn{1}{c|}{9.4}   & \multicolumn{1}{c|}{*}    \\   \hline





\hline
\end{tabular}
\caption{\noindent Performance comparison between SGDPA (two different choices of $\tau$), LALM, PDSG and FICO  in terms of CPU time (in sec.) on QCQPs having various dimensions $(n, m)$ and for strongly convex objective (first half) and convex objective (second half).}
\label{table:SGDPA}
\end{table}

\medskip 

\noindent From the preliminary numerical experiments on QCQPs with synthetic data matrices we can conclude that our algorithm is a better and viable alternative to some existing state-of-the-art methods and software. The performance of SGDPA on real data is given in the next section.

\subsection{Model predictive control  with ellipsoidal state constraints}
\label{ch07:SSS}
\noindent A linear discrete-time state-space system is a mathematical model used to describe the dynamics of a system in discrete time intervals. This model is widely used in control theory. Following the approach outlined in \cite{NedNec:14, JiaStu:05}, we examine linear discrete-time state-space systems governed by the equations:
\begin{align}\label{linearMPC}
    x(t + 1) = Ax(t) + Bu(t),
\end{align}
\noindent where \(x(t) \in \mathbb{R}^n\) is the state vector, \(u(t) \in \mathbb{R}^m\) is the control input, and \(A \in \mathbb{R}^{n\times n}\), \(B \in \mathbb{R}^{n\times m}\) are the  system matrices. The system \eqref{linearMPC} is subject to the state and input constraints:
\[
    x(t) \in \mathbb{X}_t \subseteq \mathbb{R}^n, \quad u(t) \in \mathbb{U} \subseteq \mathbb{R}^m \quad \forall t \geq 0,
\]
\noindent where we assume that the feasible set of control inputs \( \mathbb{U} \) is a simple polyhedral set  (e.g., box constraints) and the feasible set of states \( \mathbb{X}_t = \mathcal{Z}(P_t,c_t) =\{x : \|x - c_t\|^2_{P_t} \leq 1\}\) is an ellipsoidal set with the shape described by the positive semidefinite matrix $P_t$ and center $c_t$ \cite{JiaStu:05}, both sets containing the origin in their interiors. To regulate system \eqref{linearMPC} to the origin, we solve the following optimal control problem on a receding horizon of length $N$, given the current state \( x_0 \):
\begin{align}\label{eq:linearMPC2}
    \min_{\mathbf{x}, \mathbf{u}} \, & \sum_{k=0}^{N-1} \ell(x(k), u(k))  + \frac{1}{2} x^T(N) Q x(N)\\
    \text{s.t.} \quad & x(k + 1) = Ax(k) + Bu(k), \quad k = 0, \ldots, N - 1, \notag \\
    & x(k) \in \mathbb{X}_k \subset \mathbb{R}^n, \;  u(k) \in \mathbb{U} \subset \mathbb{R}^m, \quad k = 1, \ldots, N - 1,  \\
    & x(0) = x_0, \;\;  x(N) \in \mathbb{X}_N,\notag
\end{align}
where \(x(0)\) is the initial state, \(x(N)\) is the state at the end of the horizon, \( \mathbf{x} = \left(x^T(1), \ldots, x^T(N)\right)^T \), \( \mathbf{u} = \left(u^T(0), \ldots, u^T(N-1)\right)^T \), \( \ell(x, u) = \frac{1}{2} x^\top Q x + \frac{1}{2} u^\top R u \) with control cost matrix \( R \in \mathbb{R}^{m \times m} \), \( R \succeq 0 \) and state cost matrix \( Q \in \mathbb{R}^{n \times n} \), \( Q \succeq 0 \). 
The terminal set \( \mathbb{X}_N \, (:= \mathcal{Z}(P_N,c_N)) \subseteq  \mathbb{R}^n \) is also assumed to be ellipsoidal, and to contain the origin in its interior and it is added in order to ensure stability and recursive feasibility of the MPC scheme \cite{NedNec:14, JiaStu:05}.  We can eliminate the state variables by expressing the state at time \( k \) in terms of the initial state \( x_0 \) and the control inputs \( u(0), \ldots, u({k-1}) \):
\begin{equation*}
x(k) = A^k x_0 + \sum_{i=0}^{k-1} A^{k-i-1} B u(i) = A_k x_0 + B_k \mathbf{u} \quad \forall k=1:N,
\end{equation*}
where $A_k = A^k$ and $B_k = [ A^{k-1}B \ldots AB\;\; B\;\; \mathbf{0} \ldots \mathbf{0}]$.
Substituting this, the optimal control problem \eqref{eq:linearMPC2} can be restated in a compact form as:
\begin{align}\label{eq:linearMPC}
\begin{aligned}
    \min_{\mathbf{u} \in \mathbb{U}^N} \, & F (x, \mathbf{u}) \;\; \left(:= \sum_{k=0}^{N   -1} \frac{1}{2} u(k)^T R u(k) + \sum_{k=1}^{N}  \frac{1}{2} (A_k x_0  + B_k \mathbf{u})^T Q 
    (A_k x_0 + B_k \mathbf{u})  \right) \\
    \text{ s.t.} \; & \left\| A_k x_0 + B_k \mathbf{u}  - c_{k} \right\|_{P_k}^2 \leq 1  \quad \forall k=1:N. 
\end{aligned}
\end{align}
In this problem, the values of the prediction horizon $N$ can be large. Several studies highlight the performance benefits that long prediction horizons achieve in terms of closed-loop stability, harmonic distortions, and switching losses \cite{ZafLeo:23}. This feature is particularly relevant in more complex systems, such as second or third-order systems with one or more resonant frequencies \cite{KarGey:19}. This problem \eqref{eq:linearMPC} is a convex QCQP that can be solved using our algorithm SGDPA. Before presenting the test cases, let us outline the standard procedure we follow to control the system \eqref{linearMPC}, i.e., the so-called receding horizon strategy. At time $t=0$, given the initial state $x(0) =x_0$, we first solve the optimization problem \eqref{eq:linearMPC} and obtain $\mathbf{u}^*\in \mathbb{U}^N$. Using the first component of this solution, denoted as $u^*(0)\in \mathbb{R}^m$,  we compute the next state $x(1)$ using the equation $x(1) = A x_0 + Bu^*(0)$. We then set $x_0=x(1)$ and repeat this process  for a finite number of steps, known as the simulation horizon, which we denote as $N_{\text{simulation}}$. Hence, we need to solve $N_{\text{simulation}}$ times the QCQP  \eqref{eq:linearMPC} for different values of $x_0, c_k$ and $P_k$. 

\medskip

\noindent In the following, we demonstrate the performance of SGDPA algorithm through numerical experiments in two real test cases:  a second order \textit{mass-spring-damper} system \cite{EasCan:20}  and a  \textit{multimachine power system} \cite{HerAnd:11}. All the simulation parameters defining the model predictive problem \eqref{eq:linearMPC} for the mass-spring-damper system are taken from \cite{EasCan:20} and for the multimachine power system are taken from \cite{HerAnd:11}, respectively. We compare our SGDPA algorithm to its deterministic counterparts, LALM algorithm \cite{Xu:21}, and to  software package FICO \cite{FICO}. In our experiments, to solve \eqref{eq:linearMPC} using SGDPA, we consider the same stopping criteria and stepsize choice as in Section \ref{ch04:S2SS5}. Similarly, for LALM.

\begin{table}[ht]
\centering
\begin{tabular}{|c|cl|l|l|}
\hline
\multirow{2}{*}{\begin{tabular}[c]{@{}c@{}}Horizon\\ $N$\end{tabular}} & \multicolumn{2}{c|}{\begin{tabular}[c]{@{}c@{}}SGDPA \\ CPU Time \end{tabular}}                                 & \multicolumn{1}{c|}{\multirow{2}{*}{\begin{tabular}[c]{@{}c@{}}LALM\\ CPU Time\end{tabular}}} & \multicolumn{1}{c|}{\multirow{2}{*}{\begin{tabular}[c]{@{}c@{}}FICO\\ CPU Time\end{tabular}}} \\ \cline{2-3}    & \multicolumn{1}{c|}{$\tau = 0$}           & \multicolumn{1}{c|}{$\tau = 10^{-2}$} & \multicolumn{1}{c|}{}       & \multicolumn{1}{c|}{}          \\ \hline
100     & \multicolumn{1}{c|}{\begin{tabular}[c]{@{}c@{}}0.04\\ 0.03\\ 0.02\end{tabular}} &   \multicolumn{1}{c|}{\begin{tabular}[c]{@{}c@{}}\textbf{0.03}\\ \textbf{0.02}\\ \textbf{0.02 }\end{tabular}}      &    \multicolumn{1}{c|}{\begin{tabular}[c]{@{}c@{}}0.12\\ 0.1\\ 0.1\end{tabular}}      &   \multicolumn{1}{c|}{\begin{tabular}[c]{@{}c@{}}0.61\\ 0.47\\ 0.45\end{tabular}}   \\ \hline
500                                                                 & \multicolumn{1}{c|}{\begin{tabular}[c]{@{}c@{}}1.39\\ 1.19\\ 0.87\end{tabular}} &   \multicolumn{1}{c|}{\begin{tabular}[c]{@{}c@{}}\textbf{1.38}\\ \textbf{1.18}\\ \textbf{0.86}\end{tabular}}        &   \multicolumn{1}{c|}{\begin{tabular}[c]{@{}c@{}}23.2\\ 19.42\\ 15.60\end{tabular}}   &  \multicolumn{1}{c|}{\begin{tabular}[c]{@{}c@{}}8.04\\ 7.82\\ 7.64\end{tabular}} \\ \hline
1000                                                                 & \multicolumn{1}{c|}{\begin{tabular}[c]{@{}c@{}}21.86\\ 21.50\\ 20.52\end{tabular}} &   \multicolumn{1}{c|}{\begin{tabular}[c]{@{}c@{}}\textbf{21.95}\\ \textbf{21.40}\\ \textbf{20.51}\end{tabular}}        &   \multicolumn{1}{c|}{\begin{tabular}[c]{@{}c@{}}727.3\\ 716.58\\ 691.33\end{tabular}}   &  \multicolumn{1}{c|}{\begin{tabular}[c]{@{}c@{}}71.35\\ 69.15\\ 67.6\end{tabular}} \\ \hline
\end{tabular}
\caption{Mass-spring-dumper system: performance comparison between SGDPA (two choices of $\tau$), LALM and FICO in terms of CPU time (maximum time, average time and minimum time in sec.) to solve  $N_{\text{simulation}}$  QCQPs  of the form  \eqref{eq:linearMPC} for different choices of $N$. Here,  simulation horizon $N_{\text{simulation}} = 80$.}
\label{tab:Mass_spring}
\end{table}

\medskip 

\noindent Table 2 presents the  CPU times (in seconds) for the algorithms SGDPA (with two different values of $\tau$), LALM and FICO, applied to solve $N_{\text{simulation}} = 80$ QCQPs of the form \eqref{eq:linearMPC} for the mass-spring-dumper system and $N_{\text{simulation}} = 200$ QCQPs  \eqref{eq:linearMPC} for the multimachine power system, respectively. The first column provides the size of the  prediction horizons $N$, while the subsequent columns detail the computational performance of each algorithm. Specifically, each entry reports the maximum, average, and minimum CPU times for SGDPA, LALM, and FICO, respectively. As clearly illustrated, SGDPA demonstrates superior computational efficiency, consistently outperforming both LALM and FICO. This performance highlights SGDPA as a highly promising and competitive approach for solving such optimization problems (particularly, when the problems have a large number of constraints). 

\begin{table}[ht]
\centering
\begin{tabular}{|c|cl|l|l|}
\hline
\multirow{2}{*}{\begin{tabular}[c]{@{}c@{}}Horizon\\ $N$\end{tabular}} & \multicolumn{2}{c|}{\begin{tabular}[c]{@{}c@{}}SGDPA \\ CPU Time \end{tabular}}                                 & \multicolumn{1}{c|}{\multirow{2}{*}{\begin{tabular}[c]{@{}c@{}}LALM\\ CPU Time\end{tabular}}} & \multicolumn{1}{c|}{\multirow{2}{*}{\begin{tabular}[c]{@{}c@{}}FICO\\ CPU Time\end{tabular}}} \\ \cline{2-3}    & \multicolumn{1}{c|}{$\tau = 0$}           & \multicolumn{1}{c|}{$\tau = 10^{-2}$} & \multicolumn{1}{c|}{}       & \multicolumn{1}{c|}{}          \\ \hline
100     & \multicolumn{1}{c|}{\begin{tabular}[c]{@{}c@{}} $1.28$\\ $0.77$\\ $0.57$\end{tabular}} &   \multicolumn{1}{c|}{\begin{tabular}[c]{@{}c@{}}{$1.3$}\\ {$0.78$}\\ 0.58\end{tabular}}      &    \multicolumn{1}{c|}{\begin{tabular}[c]{@{}c@{}}4.71\\ 2.79\\ 2.1\end{tabular}}      &   \multicolumn{1}{c|}{\begin{tabular}[c]{@{}c@{}}\textbf{0.53}\\ \textbf{0.5}\\ \textbf{0.47}\end{tabular}}   \\ \hline
500                                                                 & \multicolumn{1}{c|}{\begin{tabular}[c]{@{}c@{}}5.49\\ 5.15\\ 4.94\end{tabular}} &   \multicolumn{1}{c|}{\begin{tabular}[c]{@{}c@{}}\textbf{5.38}\\ \textbf{5.08}\\ \textbf{4.86}\end{tabular}}        &   \multicolumn{1}{c|}{\begin{tabular}[c]{@{}c@{}} 113.76 \\ 109.22\\ 99.03\end{tabular}}   &  \multicolumn{1}{c|}{\begin{tabular}[c]{@{}c@{}}12.04\\ 10.82\\ 9.64\end{tabular}} \\ \hline
1000                                                                 & \multicolumn{1}{c|}{\begin{tabular}[c]{@{}c@{}} 34.58\\ 32.86 \\ 29.25 \end{tabular}} &   \multicolumn{1}{c|}{\begin{tabular}[c]{@{}c@{}}\textbf{33.47}\\ \textbf{30.96}\\ \textbf{27.69}\end{tabular}}        &   \multicolumn{1}{c|}{\begin{tabular}[c]{@{}c@{}} $1.08\times 10^{3}$\\ $1.02\times 10^{3}$   \\  986.21 \end{tabular}}   &  \multicolumn{1}{c|}{\begin{tabular}[c]{@{}c@{}} 97.71\\ 95.28\\ 93.81  \end{tabular}} \\ \hline
\end{tabular}
\caption{Multimachine power system: performance comparison between SGDPA (two choices of $\tau$), LALM and FICO in terms of CPU time (maximum time, average time and minimum time in sec.) to solve  $N_{\text{simulation}}$  QCQPs  of the form  \eqref{eq:linearMPC} for different choices of $N$. Here,  simulation horizon $N_{\text{simulation}} = 200$.}
\label{tab:Power_system}
\end{table}




\section{Conclusions}
In this paper we have addressed the design and the convergence analysis of an  optimization algorithm for solving efficiently smooth convex problems with a large number of constraints. To achieve these goals, we have proposed a stochastic perturbed augmented Lagrangian (SGDPA) method, where a perturbation is introduced in the augmented Lagrangian function by multiplying the dual variables with a subunitary parameter. Due to the computational simplicity of the SGDPA subproblem, which uses only gradient information and considers only one constraint at a time, makes our algorithm suitable for problems with many functional constraints. We have provided a detailed convergence analysis in both optimality and feasibility criteria for the iterates of SGDPA algorithm using basic assumptions on the problem. Preliminary numerical experiments on  problems with many quadratic constraints have demonstrated the viability and performance of our method when compared to some existing state-of-the-art optimization methods and software.



\end{document}